\documentclass[a4paper,oneside,12pt]{article}
\usepackage[utf8]{inputenc}
\usepackage{lmodern}
\usepackage[T1]{fontenc}
\usepackage{verbatim}
\usepackage{textcomp}
\usepackage[english]{babel}
\usepackage[a4paper,vmargin={2.5cm,2.5cm},hmargin={1.5cm,1.5cm}]{geometry}
\usepackage[font=sf, labelfont={sf,bf}, margin=1cm]{caption}
\usepackage[pdftex]{hyperref}
\usepackage[pdftex]{color,graphicx}
\usepackage{amssymb}
\usepackage{amsmath}
\def\llbracket{[\hspace{-.10em} [ }

\def\rrbracket{ ] \hspace{-.10em}]}
\renewcommand{\leq}{\leqslant}
\renewcommand{\geq}{\geqslant}

\usepackage{amsmath,amsfonts,amssymb,amsthm,mathrsfs}

\def\build#1_#2^#3{\mathrel{
\mathop{\kern 0pt#1}\limits_{#2}^{#3}}}

\renewcommand{\d}{\mathrm{d}}
\newcommand{\ind}{\mathbf{1}}
\newcommand{\E}{\mathbb{E}}
\theoremstyle{plain}
\newtheorem{theorem}{Theorem}
\newtheorem*{thm}{Theorem} 
\newtheorem{corollary}{Corollary}
\newtheorem{proposition}[corollary]{Proposition}

\newtheorem{lemma}[corollary]{Lemma}
\newtheorem{open}[corollary]{Question}
\theoremstyle{definition}
\newtheorem{definition}[corollary]{Definition}

\theoremstyle{remark}
\newtheorem*{remark}{Remark}
\newcommand{\op}[1]{\operatorname{#1 }}
\renewcommand{\P}[1]{\mathbb{P}\left(#1\right)}

  \author{N.\ Curien, L.\ Ménard, G.\ Miermont} \title{\bf A view from
    infinity of the uniform infinite planar quadrangulation}
  \date{{\small \today}}

\newcommand*\ltrees[2]{\mathbf{T}_{#1}^{(#2)}}

\newcommand*\wltrees[2]{\overline{\mathbf{T}} \ \hspace{-1.2mm}^{(#2)}_{#1}}

\newcommand*\carte{\mathbf{Q}}
\newcommand*\arbres{\mathbf{T}}

\begin{document}

\maketitle

\abstract{We introduce a new construction of the Uniform Infinite Planar Quadrangulation (UIPQ). Our approach is based on an extension of the Cori-Vauquelin-Schaeffer mapping in the context of infinite trees, in the spirit of previous work \cite{CD06,LGM10,Men08}. However, we release the positivity constraint on the labels of trees which was imposed in these references, so that our construction is technically much simpler. This approach allows us to prove the conjectures of Krikun \cite{Kri08} pertaining to the ``geometry at infinity'' of the UIPQ, and to derive new results about the UIPQ, among which  a fine study of infinite geodesics.}


\section{Introduction}
\label{sec:introduction}

The purpose of this work is to develop a new approach to the Uniform Infinite Planar Quadrangulation (UIPQ), which is a model of random discrete planar geometry consisting in a cell decomposition of the plane into quadrangles, chosen ``uniformly at random'' among all homeomorphically distinct possibilities.

Recall that a planar map is a proper embedding of a finite connected graph in the two-dimensional sphere, viewed up to orientation-preserving homeomorphisms of the sphere.  The faces are the connected components of the complement of the union of the edges. A map is a triangulation (respectively a quadrangulation) if all its faces are incident to three (respectively four) edges. A map is rooted if one has distinguished an oriented edge called the root edge. Planar maps are basic objects in combinatorics and have been extensively studied since the work of Tutte in the sixties \cite{Tut63}. They also appear in various areas, such as algebraic geometry \cite{LZ04}, random matrices \cite{Zvo97} and theoretical physics, where they have been used as a model of random geometry \cite{ADJ97}. 

The latter was part of Angel and Schramm's motivation to introduce in \cite{AS03} the so-called Uniform Infinite Planar Triangulation as a model for random planar geometry.  Its companion, the UIPQ, was later defined by Krikun \cite{Kri05} following a similar approach. One advantage of quadrangulations over triangulations is that there exists a very nice bijection between, on the one hand, rooted planar quadrangulations with $n$ faces, and on the other hand, labeled plane trees with $n$ edges and non-negative labels. This bijection is due to Cori and Vauquelin \cite{CV81}, but only reached its full extension with the work of Schaeffer \cite{Sch98}. See Section \ref{sec:schaeffer}. This leads Chassaing and Durhuus \cite{CD06} to introduce an infinite random quadrangulation of the plane, generalizing the Cori-Vauquelin-Schaeffer bijection to a construction of the random quadrangulation from an infinite random tree with non-negative labels. Ménard \cite{Men08} then showed that the two constructions of \cite{Kri05,CD06} lead to the same random object. In the Chassaing-Durhuus approach, the labels in the random infinite tree correspond to distances from the origin of the root edge in the quadrangulation, and thus information about the labels can be used to derive geometric properties such as volume growth around the root in the UIPQ \cite{CD06,LGM10,Men08}.

Let us describe quickly the UIPQ with the point of view of Angel-Schramm and Krikun.  
If $Q_n$ is a random rooted quadrangulation uniformly distributed over the set of all rooted quadrangulations with $n$ faces, then we have \cite{Kri05}
\begin{eqnarray*} Q_n &\xrightarrow[n\to \infty]{(d)}& Q_\infty, \end{eqnarray*} in distribution in the sense of the {\em local convergence}, meaning that for every fixed $R>0$, the combinatorial balls of $Q_n$ with radius $R$ and centered at the root converge in distribution as $n\to\infty$ to that of $Q_\infty$, see Section \ref{sec:uipq} for more details.
The object $Q_\infty$ is a random infinite rooted quadrangulation called the Uniform Infinite Planar Quadrangulation (UIPQ).  The UIPQ and its sister the UIPT are fundamental objects in random geometry and have been the object of many studies. See \cite{Ang03,AS03,Ben10,Kri05,Kri08,Kri04} and references therein.

In the present work, we give a new construction of the UIPQ from a certain random labeled tree. This is in the spirit of the ``bijective'' approach by Chassaing-Durhuus, but where the positivity constraint on the labels is released. Though the labels no longer correspond to distances from the root of the UIPQ, they can still be interpreted as ``distances seen from the point at infinity''. In many respects, this construction is simpler than \cite{CD06} because the unconditioned labeled tree has a very simple branching structure --- its genealogy is that of a critical branching process conditioned on non-extinction. This simplifies certain computations on the UIPQ and enables us to derive new results easily.


Let us briefly describe our construction.%
 We denote by $T_{\infty}$ the critical geometric Galton-Watson tree conditioned to survive. This random infinite planar tree with one end has been introduced by Kesten \cite{Kes86} and can be built from a semi-infinite line of vertices $x_0,x_1,x_2,\ldots$ together with independent critical geometric Galton-Watson trees grafted to the left-hand side and right-hand side of each vertex $x_i$ for  $i \geq 0$, see Section \ref{sec:uiltuipq}. Conditionally on $T_{\infty}$, we consider a sequence of independent variables $(\mathsf{d}_{e})_{e \in E(T_{\infty})}$ indexed by the edges of $T_{\infty}$ which are uniformly distributed over $\{-1,0,+1\}$. We then assign to every vertex $u$ of $T_{\infty}$ a label $\ell(u)$ corresponding to the sum of the numbers $ \mathsf{d}_{e}$ along the ancestral path from $u$ to the root $x_{0}$ of $T_{\infty}$. Given an extra Bernoulli variable $\eta \in \{0,1\}$ independent of $(T_{\infty},\ell)$, it is then possible to extend the classical Schaeffer construction to define a quadrangulation $\Phi((T_{\infty},\ell),\eta)$ from $(T_{\infty},\ell)$ and $\eta$, see Section \ref{sec:schaeffer}. The only role of $\eta$ is to prescribe the orientation of the root edge in $\Phi((T_{\infty},\ell),\eta)$. The random infinite rooted quadrangulation $Q_\infty=\Phi((T_{\infty},\ell),\eta)$ has the distribution of the UIPQ, see Theorem \ref{sec:unif-infin-label-1}. Moreover, the vertices of $Q_{\infty}$ correspond to those of $T_{\infty}$ and via this identification,  Theorem \ref{sec:unif-infin-label-1} gives a simple interpretation of the labels: Almost surely, for any pair of vertices $u,v$ of $Q_{\infty}$ we have
 \begin{eqnarray*} 
   \ell(u)-\ell(v) &=& \lim_{z \to \infty} \left(
\op{d}_{\op{gr}}^{Q_{\infty}}(u,z) -
\op{d}_{\op{gr}}^{Q_{\infty}}(v,z) \right)\, , \quad \quad (*)
 \end{eqnarray*}
where $\op{d}_{\op{gr}}^{Q_{\infty}}$ is the usual graph distance. 
 The fact that the limit exists as $z\to \infty$ in $(*)$ means that the right-hand side is constant everywhere but on a finite subset of vertices of $Q_{\infty}$.  Theorem \ref{sec:unif-infin-label-1} and its corollaries also answer positively the three conjectures raised by Krikun in \cite{Kri08}. Note that the existence of the limit in $(*)$ was shown  in \cite{Kri08}. It also follows from our fine study of the geodesics and their coalescence properties in the UIPQ, see Proposition \ref{sec:confluent-geodesics-1} and Theorem \ref{cut-points}.

 As a corollary of our new construction we 
 study (see Theorem \ref{separatingcycles}) the length of the separating cycle at a given height (seen from $\infty$) in the UIPQ, much in the spirit of a previous work of Krikun's \cite{Kri05}. We also deduce new properties that support a conjecture of Angel \& Schramm \cite{AS03} (reformulated in our context) saying that the UIPQ is recurrent. Namely, we show that the distances from infinity along the random walk on the UIPQ is a recurrent process.
\bigskip
 
 The paper is organized as follows. In Section $2$ we introduce the construction of the UIPQ based on a random infinite labeled tree and present our main theorem. Section $3$ is devoted to the proof of Theorem \ref{sec:unif-infin-label-1}, which goes through an analysis of discrete geodesics in the UIPQ. In particular, we establish  a confluence property of geodesics towards the root (Proposition \ref{sec:confluent-geodesics-1}) and  a certain uniqueness property of geodesic rays towards infinity (Theorem \ref{cut-points}). Section \ref{scaling-limits} is devoted to the study of the scaling limits for the contour functions describing the infinite labeled tree $(T_{\infty},\ell)$ and to the proofs of two technical lemmas used to derive Theorem \ref{cut-points}. Using our new construction we finally study separating cycles at a given heigh (Section \ref{sec:horoballs}) and random walk on the UIPQ (Section \ref{SRW}). 
 
\bigskip
 
\noindent \textbf{Acknowledgments:} We deeply thank Jean-Fran\c cois Le Gall for fruitful discussions and a careful reading of a first version of this article.

\section{The UIPQ and the uniform infinite labeled tree}\label{sec:uipq-its-encoding}

\subsection{Finite and infinite quadrangulations}
\label{subsec:quadrangulations}

Consider a proper embedding of a finite connected graph in the sphere
$\mathbb{S}_2$ (loops and multiple edges are allowed).  A \emph{finite
  planar map} $m$ is an equivalence class of such embeddings modulo orientation preserving homeomorphisms of the sphere. Let $\overrightarrow{E}(m)$ be the set of all oriented edges of $m$ (each edge corresponds to exactly two oriented edges). A planar map is \emph{rooted} if it has a distinguished oriented edge $e^* \in \overrightarrow{E}(m)$, which is called the root edge.
 If $e$ is an oriented edge of a map we write $e_-$ and $e_+$ for its origin and target vertices and $\overleftarrow{e}$ for the reversed edge.

The set of vertices of a map $m$ is denoted by $V(m)$. We will equip $V(m)$ with the graph distance: If $v$ and $v'$ are two vertices, $\op{d}_{\op{gr}}^m(v,v')$ is the
minimal number of edges on a path from $v$ to $v'$ in $m$. If $v\in V(m)$, the \emph{degree} of $v$ is the number of oriented edges pointing towards $v$ and is denoted by $\op{deg}(v)$.
 
 The {\em
   faces} of the map are the connected components of the complement of the union of its edges. The {\em degree} of a face is the number of edges that are incident to it, where it should be understood that an edge lying entirely in a face is incident twice to this face.  A finite planar map is a \emph{quadrangulation} if all its faces have degree $4$, that is $4$ incident edges. A planar map is a \emph{quadrangulation with holes} if all its faces have degree $4$, except for a number of distinguished faces which can be of arbitrary even degrees. We call these faces the {\em holes}, or the {\em boundaries} of the quadrangulation. 

\subsubsection{Infinite quadrangulations and their planar embeddings}\label{sec:infin-quadr}

Let us introduce infinite quadrangulations using the approach of
Krikun \cite{Kri05}, see also \cite{AS03,BS01}. For every integer $n \geq 1$, we denote by
$\carte_n$ the set of all rooted quadrangulations with $n$ faces.  For
every pair $q,q' \in \carte_f = \bigcup_{n \geq 1} \carte_n$ we define
 \begin{eqnarray*} 
d_{\carte} \left( q,q' \right) &=& \Big( 1 + \sup \left\{ r: \,
B_{\carte,r}(q) = B_{\carte,r}(q') \right\} \Big)^{-1}
 \end{eqnarray*}
where, for $r \geq 1$, $B_{\carte,r}(q)$ is the planar map whose edges (resp. vertices) are all edges (resp. vertices) incident to a face of $q$ having at least one vertex at
distance strictly smaller than $r$ from the root vertex $e^*_{-}$, and
$\sup \emptyset = 0$ by convention. Note that $B_{\carte,r}(q)$ is a
quadrangulation with holes. 

The pair $(\carte_f,d_{\carte})$ is a metric space, we let
$(\carte,d_{\carte})$ be the completion of this space. We call
\emph{infinite quadrangulations} the elements of $\carte$ that are not
finite quadrangulations and we denote the set of all such
quadrangulations by $\carte_{\infty}$.  Note that one can extend the
function $q \in \carte_f \mapsto B_{\carte,r}(q)$ to a continuous
function $B_{\carte,r}$ on $\carte$.

\paragraph{Infinite quadrangulations of the plane.} An infinite quadrangulation $q$ defines a unique infinite graph $G$ with a
root edge, together with a consistent family of planar embeddings
$(B_{\carte,r}(q),r\geq 1)$ of the combinatorial balls of $G$ centered
at the root vertex. 

Conversely, any sequence $q_1,q_2,\ldots$ of rooted quadrangulations
with holes, such that $q_r=B_{\carte,r}(q_{r+1})$ for every
$r\geq 1$, specifies a unique infinite quadrangulation $q$ whose ball
of radius $r$ is $q_r$ for every $r\geq 1$.
\begin{definition} An infinite quadrangulation $q \in \carte_{\infty}$ is called a \emph{quadrangulation of the plane} if it has one end, that is, if for any $r \geq 0$ the graph $q \backslash B_{\carte,r}(q)$ has only one infinite connected component.
\end{definition}
It is not hard to convince oneself that quadrangulations of the plane also coincide with equivalence classes of certain proper embeddings of infinite graphs in the plane $\mathbb{R}^2$, viewed up to orientation preserving homeomorphisms. Namely these are the proper embeddings $\chi$ of locally finite planar graphs such that
\begin{itemize}
\item every compact subset of $\mathbb{R}^2$ intersects only finitely many edges of $\chi$,
 \item the connected components of the complement of the union of edges of $\chi$ in $\mathbb{R}^2$ are all bounded topological quadrangles.
 \end{itemize}
 \begin{remark} Note that a generic element of $\carte_{\infty}$ is not necessarily a quadrangulation of the plane. See \cite{AS03,CD06,Men08} and the Appendix below for more details about this question. 
 \end{remark}

\subsubsection{The Uniform Infinite Planar Quadrangulation}\label{sec:uipq}

Now, let $Q_n$ be a random variable with uniform distribution on
$\carte_n$. Then as $n\to\infty$, the sequence $(Q_n)_{n\geq 1}$ converges in
distribution to a random variable with values in $\carte_\infty$. 
\begin{thm}[{\cite{Kri05}}]\label{krikun} 
For every $n \geq 1$, let $\nu_n$ be the uniform probability measure on
$\mathbf{Q}_n$. The sequence $(\nu_n)_{n \geq 1}$ converges to
a probability measure $\nu$, in the sense of weak convergence in the
space of probability measures on $\left(\mathbf{Q}, d_{\mathbf{Q}}
\right)$. Moreover, $\nu$ is supported on the set of infinite rooted
quadrangulations of the plane.
\end{thm}
 The probability measure $\nu$ is called the law of the uniform
infinite planar quadrangulation (UIPQ). 

\subsection{Labeled trees}\label{sec:spatial-trees}

Throughout this work we will use the standard formalism for planar
trees as found in \cite{Nev86}. Let
 \begin{eqnarray*}
\mathcal{U}& = &\bigcup_{n=0}^{\infty} \mathbb{N}^n
 \end{eqnarray*}
where $\mathbb{N} = \{ 1,2, \ldots \}$ and $\mathbb{N}^0 = \{
\varnothing \}$ by convention. An element $u$ of $\mathcal{U}$ is thus
a finite sequence of positive integers. If $u, v \in \mathcal{U}$,
$uv$ denotes the concatenation of $u$ and $v$. If $v$ is of the form
$uj$ with $j \in \mathbb{N}$, we say that $u$ is the \emph{parent} of
$v$ or that $v$ is a \emph{child} of $u$. More generally, if $v$ is of
the form $uw$, for $u,w \in \mathcal{U}$, we say that $u$ is an
\emph{ancestor} of $v$ or that $v$ is a \emph{descendant} of $u$. A
\emph{rooted planar tree} $\tau$ is a (finite or infinite) subset of
$\mathcal{U}$ such that
\begin{enumerate}
\item $\varnothing \in \tau$ ($\varnothing$ is called the \emph{root}
  of $\tau$),
\item if $v \in \tau$ and $v \neq \varnothing$, the parent of $v$
  belongs to $\tau$
\item for every $u \in \mathcal{U}$ there exists $k_u(\tau) \geq 0$
  such that $uj \in \tau$ if and only if $j \leq k_u(\tau)$.
\end{enumerate}
A rooted planar tree can be seen as a graph, in which an edge links
two vertices $u,v$ such that $u$ is the parent of $v$ or vice-versa.
This graph is of course a tree in the graph-theoretic sense, and has a
natural embedding in the plane, in which the edges from a vertex $u$
to its children $u1,\ldots,uk_u(\tau)$ are drawn from left to right.

We let $|u|$ be the length of the word $u$. 
The
integer $|\tau|$ denotes the number of edges of $\tau$ and is called
the size of $\tau$.  
A \emph{spine} in a tree $\tau$ is an infinite sequence
$u_{0},u_{1},u_{2},\ldots$ in $\tau$ such that $u_{0}=\varnothing$ and $u_i$
is the parent of $u_{i+1}$ for every $i\geq 0$. If $a$ and $b$ are two
vertices of a tree $\tau$, we denote the set of vertices along
the unique geodesic path going from $a$ to $b$ in $\tau$ by $\llbracket a,b \rrbracket$.

\bigskip

A \emph{rooted labeled tree} (or spatial tree) is a pair $\theta =
(\tau, (\ell(u))_{u \in \tau})$ that consists of a rooted planar tree $\tau$
and a collection of integer labels assigned to the vertices of $\tau$,
such that if $u,v \in \tau$ and $v$ is a child of $u$, then $|\ell(u)
- \ell(v)| \leq 1$.  For every $l \in \mathbb{Z}$, we denote by
$\ltrees{}{l}$ the set of labeled trees for which $\ell(\varnothing) =
l$, and $\ltrees{\infty}{l}$, resp. $\ltrees{f}{l}$, resp. $\ltrees{n}{l}$, are the subsets of $\ltrees{}{l}$ consisting of the infinite trees,  resp. finite trees, resp. trees with $n$ edges.  If $\theta =
(\tau,\ell)$ is a labeled tree, $|\theta| = |\tau|$ is the size of
$\theta$.

\bigskip

If $\theta$ is a labeled tree and $h \geq 0$ is an integer, we denote the labeled subtree of $\theta$ consisting of all vertices of $\theta$ and their labels up to height $h$ by $B_{\arbres,h}(\theta)$. For every pair $\theta, \theta'$ of labeled trees define
 \begin{eqnarray*}
d_{\arbres}(\theta,\theta') &= &\big( 1 + \sup \left\{ h: \,
B_{\arbres,h}(\theta) = B_{\arbres,h}(\theta') \right\} \big)^{-1}.
 \end{eqnarray*}
One easily checks that $d_{\arbres}$ is a distance on  the set of all labeled trees,  which turns this set into a separable and complete metric space. 

In the rest of this work we will mostly be  interested in the following set of infinite trees. We let
$\mathscr{S}$ be the set of all labeled trees $(\tau,\ell)$ in
$\ltrees{\infty}{0}$ such that
\begin{itemize}
\item $\tau$ has exactly one spine, which we denote by
  $\varnothing = \mathrm{S}_{\tau}(0), \mathrm{S}_{\tau}(1), \mathrm{S}_{\tau}(2), \ldots $
\item
$\inf_{i\geq 0}\ell(\mathrm{S}_{\tau}(i))=-\infty$.
\end{itemize}
If $\theta=(\tau,\ell) \in \mathscr{S}$, the spine then splits $\tau$ in two parts, which we call the left and right parts, and every vertex $\mathrm{S}_\tau(n)$ of the spine determines a subtree of $\tau$ to its left and to its right. These are denoted by $L_n(\theta),R_n(\theta)$, formally, 
\begin{eqnarray*}
L_n(\theta)&=&\{v\in \mathcal{U}:\mathrm{S}_{\tau}(n)v\in \tau,\mathrm{S}_{\tau}(n)v\prec
\mathrm{S}_{\tau}(n+1)\}\\ R_n(\theta)&=& \{v\in \mathcal{U}:\mathrm{S}_{\tau}(n)v\in \tau,
\mathrm{S}_{\tau}(n+1)\prec \mathrm{S}_{\tau}(n)v\}\cup \{\varnothing\}\, ,
\end{eqnarray*}
where $u\prec v$ denotes the lexicographical order on $\mathcal{U}$.  The subtrees $L_n(\theta),R_n(\theta)$ naturally inherit the labels from $\theta$, so that we really see $L_n(\theta),R_n(\theta)$ as elements of $\ltrees{f}{X_n(\theta)}$, where $X_n(\theta)=\ell(\mathrm{S}_\tau(n))$ is the label of the $n$-th vertex of the spine.  We can of course reconstruct the tree $\theta$ from the sequence $(L_n(\theta),R_n(\theta))_{n\geq 0}$.  In the sequel, we will often write $ \mathrm{S}(n),X_n,L_n,R_n$ instead of $ \mathrm{S}_{\tau}(n),X_n(\theta), L_n(\theta),R_n(\theta)$ when there is no ambiguity on the underlying labeled tree.


\subsection{The Schaeffer correspondence}
\label{sec:schaeffer}
One of the main tools for studying random quadrangulations is a bijection initially due to Cori \& Vauquelin \cite{CV81}, and that was much developed by Schaeffer \cite{Sch98}. It establishes a one-to-one correspondence between rooted and pointed quadrangulations with $n$ faces, and pairs consisting of a labeled tree of $\mathbf{T}_{n}^{(0)}$ and an element of $\{0,1\}$. Let us describe this correspondence and its extension to infinite quadrangulations.

\subsubsection{From trees to quadrangulations}\label{sec:from-trees-quadr}
A rooted and \emph{pointed} quadrangulation is a pair 
$\mathbf{q}=(q,\rho)$ where $q$ is a rooted quadrangulation and $\rho$
is a distinguished vertex of $q$. We write $\carte_n^{\bullet}$ for
the set of all rooted and pointed quadrangulations with $n$ faces. We first describe the mapping from labeled trees to quadrangulations.

Let $\theta=(\tau,\ell)$ be an element of $\ltrees{n}{0}$. We view
$\tau$ as embedded in the plane. A \emph{corner} of a vertex in $\tau$
is an angular sector formed by two consecutive edges in clockwise
order around this vertex. Note that a vertex of degree $k$ in $\tau$
has exactly $k$ corners. If $c$ is a corner of $\tau$, $\mathcal{V}(c)$
denotes the vertex incident to $c$. By extension, the label $\ell(c)$ of a
corner $c$ is the label of $\mathcal{V}(c)$.

The corners are ordered clockwise cyclically around the tree in the so-called {\em contour order}. If we view $\tau$ as a planar map with one face, then $\tau$ can be seen as a polygon with $2n$ edges that are glued by pairs, and the contour order is just the usual cyclic order of the corners of this polygon. 
We fix the labeling by letting $(c_{0},c_{1},c_{2}, \ldots,c_{2n-1})$ be the sequence of
corners visited during the contour process of $\tau$, starting from
the corner $c_0$ incident to $\varnothing$ that is located to the left
of the oriented edge going from $\varnothing$ to $1$ in $\tau$. We extend this sequence of corners into a sequence $(c_i,i\geq
0)$ by periodicity, letting $c_{i+2n}=c_i$.  For $i \in \mathbb{Z}_{+}$,
the \emph{successor} $\mathcal{S}(c_{i})$ of $c_{i}$ is the first corner $c_{j}$ in the
list $c_{i+1},c_{i+2},c_{i+3},\ldots$ of label $\ell(c_j)=\ell(c_{i})-1$, if such a corner exists. In
the opposite case, the successor of $c_i$ is an extra element
$\partial$, not in $\{c_i,i\geq 0\}$.

Finally, we construct a new graph as follows. Add an extra vertex $\rho$
in the plane, that does not belong to (the embedding of) $\tau$. For 
every corner $c$, draw an arc between $c$ and its successor if this
successor is not $\partial$, or draw an arc between $c$ and $\rho$ if
the successor of $c$ is $\partial$. The construction can be made in such a way that the arcs do not cross.  After the interior of the edges of
$\tau$ has been removed, the resulting embedded graph, with vertex set
$\tau\cup\{\rho\}$ and edges given by the newly drawn arcs, is a
quadrangulation $q$. In order to root this quadrangulation, we consider
some extra parameter $\eta\in \{0,1\}$. If $\eta=0$, the root of $q$
is the arc from $c_0$ to its successor, oriented in this direction. If
$\eta=1$ then the root of $q$ is the same edge, but with opposite
orientation.  We let $q=\Phi(\theta,\eta)\in \carte_n^\bullet$
($q$ comes naturally with the distinguished vertex $\rho$).

\begin{thm}[Theorem 4 in \cite{CS04}] 
\label{schaeffer}
The mapping $\Phi : \mathbf{T}^{(0)}_{n} \times\{0,1\}
\longrightarrow \carte_{n}^{\bullet}$ is a bijection. If
$q=\Phi((\tau,\ell),\eta)$ then for every vertex $v$ of $q$ not
equal to $\rho$, one has
 \begin{eqnarray}\label{eq:1}
\op{d}^q_{\mathrm{gr}}(v,\rho) &=&  \ell(v)-\min_{u\in \tau}\ell(u)+1\, ,
 \end{eqnarray}
  where we recall that every vertex of $q$ not equal to $\rho$ is identified to a vertex
  of $\tau$.
\end{thm}

Note that (\ref{eq:1}) can also be rewritten as 
 \begin{eqnarray}\label{eq:2}
\ell(v)&=&\op{d}^q_{\mathrm{gr}}(v,\rho)-\op{d}^q_{\mathrm{gr}}(e^*_{\pm},\rho)\,
,\qquad v\in V(q)\, ,
 \end{eqnarray}
where 
$$e^*_\pm=\varnothing=\left\{\begin{array}{cl} e^*_- &\mbox{ if
}\op{d}^q_{\mathrm{gr}}(e^*_-,\rho)-\op{d}^q_{\mathrm{gr}}(e^*_+,\rho)=-1\\ e^*_+
&\mbox{ if }\op{d}^q_{\mathrm{gr}}(e^*_-,\rho) -\op{d}^q_{\mathrm{gr}}(e^*_+,\rho)= 1
\end{array}
\right.$$ Hence, these labels can be recovered from the pointed
quadrangulation $(q,\rho)$. This is of course not surprinsing since the
function $\Phi:\ltrees{n}{0}\times\{0,1\}\to \carte_n^\bullet$ is
invertible (see the next section for the description of the inverse mapping).

\begin{figure}[h]

\begin{center} 
\includegraphics[width=15cm]{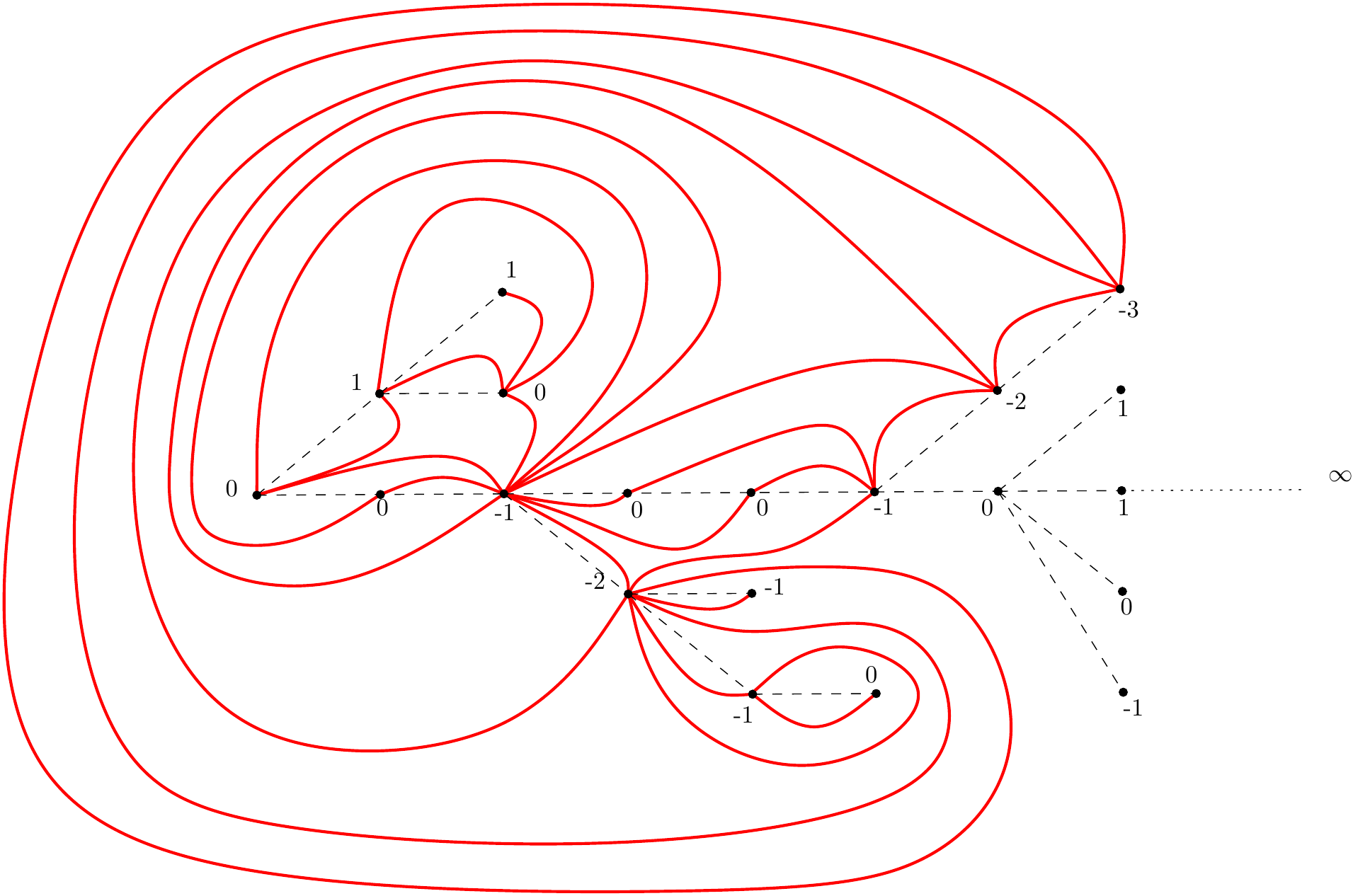}
\end{center}
\caption{ \label{fig:schaeffer2} Illustration of the Schaeffer correspondence. The tree is represented in dotted lines and  the quadrangulation in solid lines.}
\end{figure}

\paragraph{Infinite case.} 
We now aim at extending the construction of $\Phi$ to elements of $ \mathscr{S}$. Let $(\tau, (\ell(u))_{u\in\tau})\in \mathscr{S}$.  Again, we consider an embedding of $\tau$ in the plane, with isolated vertices. This is always possible (since $\tau$ is locally finite). The notion of a corner is unchanged in this setting, and there is still a notion of clockwise contour order for the corners of $\tau$, this order being now a total order, isomorphic to $(\mathbb{Z},\leq)$, rather than a cyclic order.  We consider the sequence $(c_{0}^{(L)},c_{1}^{(L)},c_{2}^{(L)}, \ldots)$ of corners visited by the contour process of the left side of the tree in clockwise order --- roughly speaking, these corners correspond to the concatenation of the contour orders around the trees $L_n(\theta),n\geq 0$, plus the extra corners induced by grafting these trees on the spine. Similarly, we denote the sequence of corners visited on the right side by $(c_{0}^{(R)},c_{1}^{(R)},c_{2}^{(R)},\ldots)$, in counterclockwise order. Notice that $c_{0}^{(L)}=c_{0}^{(R)}$ denotes the corner where the tree has been rooted. We now concatenate these two sequences into a unique sequence indexed by $\mathbb{Z}$, by letting, for $i\in \mathbb{Z}$,
$$c_{i}=\left\{\begin{array}{cc}
c_i^{(L)} & \mbox{ if }i\geq 0\\
c_{-i}^{(R)} & \mbox{ if }i< 0\, .
\end{array}\right.$$ 
In the sequel, we will write $c_i\leq c_j$ if $i\leq j$.  
For any $i \in \mathbb{Z}$, the \emph{successor} $\mathcal{S}(c_{i})$ of $c_{i}$ is the
first corner $c_{j}\geq c_{i+1}$ such
that the label $\ell(c_{j})$ is equal to
$\ell(c_{i})-1$. From the assumption that $\inf_{i\geq 0}
\ell(\mathrm{S}_{\tau}(i))=-\infty$, and since all the vertices of the spine appear in the sequence $(c_{i}^{(L)})_{i\geq 0}$,
it holds that each corner has exactly one successor. We can associate
with $(\tau,(\ell(u))_{u\in \tau})$ an embedded graph $q$ by drawing an arc
between every corner and its successor. See
Fig. \ref{fig:schaeffer2}. Note that, in contrast with the above
description of the Schaeffer bijection on $\ltrees{n}{0}\times
\{0,1\}$, we do not have to add an extra distinguished vertex $\rho$
in this context.

In a similar way as before, the embedded graph $q$ is rooted at the
edge emerging from the distinguished corner $c_{0}$ of $\varnothing$,
that is, the edge between $c_{0}$ and its successor $\mathcal{S}(c_{0})$. The
direction of the edge is given by an extra parameter $\eta \in
\{0,1\}$, similarly as above.

\begin{proposition} \label{prop:continuity} The resulting embedded graph $q$ is an infinite 
quadrangulation of the plane, and the extended mapping
$\Phi:\mathscr{S}\cup \ltrees{f}{0}\to \carte$ is continuous.
\end{proposition}

\proof We first check that every corner in $\tau$ is the successor of
only a finite set of other corners. Indeed, if $c$ is a corner,
say $c=c_i$ for $i\in \mathbb{Z}$, then from the assumption that
$\inf_j\ell(\mathrm{S}_{\tau}(j))=-\infty$, there exists a corner $c_j$ with $j<i$ such
that the vertex incident to $c_j$ belongs to the spine
$\{\mathrm{S}_{\tau}(0),\mathrm{S}_{\tau}(1),\ldots\}$, and $\min_{j\leq k\leq i}
\ell(c_k)<\ell(c_i)-1$. Therefore, for every $k\leq j$, the successor
of $c_k$ is not $c_i$. 

Together with the fact that every vertex has a number of successors equal to its degree, this shows that the embedded graph $q$ is locally finite, in the sense that every vertex is incident to a finite number of edges. The fact that every face of $q$ is a quadrangle is then a consequence of the construction of the arcs, as proved e.g.\ in \cite{CS04}. It remains to show that $q$ can be properly embedded in the plane, that is, has one end. This comes from the construction of the edges and the fact that $\tau$ has only one end. The details are left to the reader.

To prove the continuity of $\Phi$, let $\theta_n=(\tau_n,\ell_n)$
be a sequence in $\mathscr{S}\cup \ltrees{f}{0}$ converging to
$\theta=(\tau,\ell)\in \mathscr{S}\cup \ltrees{f}{0}$. If $\theta\in
\ltrees{f}{0}$ then $\theta_n=\theta$ for every $n$ large enough, so
the fact that $\Phi(\theta_n)\to \Phi(\theta)$ is obvious. So
let us assume that $\theta\in \mathscr{S}$, with spine vertices
$\mathrm{S}_{\tau}(0),\mathrm{S}_{\tau}(1),\ldots$. Let $R>0$ be an integer, and let $l(R)$ be
the minimal label of a vertex in $B_{\arbres,R}(\theta)$. Since
$\inf(\ell(\mathrm{S}_{\tau}(i)))=-\infty$, we can define $f(R)> R$ as the first $i\geq 1$
such that $\ell(\mathrm{S}_{\tau}(i))=l(R)-2$. If $c$ is a corner in the subtree of
$\tau$ above $\mathrm{S}_{\tau}(f(R))$, then the successor of $c$ cannot be in
$B_{\arbres,R}(\theta)$. Indeed, if $\ell(c)\geq l(R)-1$ then the successor of $c$ has
to be also in the subtree of $\tau$ above $\mathrm{S}_{\tau}(f(R))$, while if
$\ell(c)<l(R)-1$, then this successor also has label $<l(R)-1$, and
thus cannot be in $B_{\arbres,R}(\theta)$ by definition. Similarly, $c$ cannot
be the successor of any corner in $B_{\arbres,R}(\theta)$, as these successors
 necessarily  are in the subtree of $\tau$ below $S_{\tau}(f(R))$.

Now, for every $n$ large enough, it holds that
$B_{\arbres,f(R)}(\theta_n)=B_{\arbres,f(R)}(\theta)$, from which we obtain that the
maps formed by the arcs incident to the vertices of
$B_{\arbres,R}(\theta)=B_{\arbres,R}(\theta_n)$ are the same, and moreover, no extra arc
constructed in $\theta_n$ or $\theta$ is incident to a vertex of
$B_{\arbres,R}(\theta)=B_{\arbres,R}(\theta_n)$. Letting $r>0$ and choosing $R$ so that
 all the edges of $B_{\carte,r}(\Phi(\theta))$ appear as arcs incident to vertices of $B_{\arbres,R}(\theta)$, we obtain that
$B_{\carte,r}(\Phi(\theta))=B_{\carte,r}(\Phi(\theta_n))$ for
$n$ large enough. Therefore, we get that $\Phi(\theta_n)\to
\Phi(\theta)$, as desired.  \endproof

The vertex set of $q$ is precisely $\tau$, so that the labels $\ell$
on $\tau$ induce a labeling of the vertices of $q$. In the finite
case, we saw earlier in (\ref{eq:2}) that these labels could be
recovered from the pointed quadrangulation obtained from a finite
labeled tree. In our infinite setting, this is much less obvious: Intuitively the distinguished vertex $\rho$ of the finite case is ``lost at
infinity''.

\paragraph{Bounds on distances.}

We will see later that when the
infinite labeled tree has a special distribution corresponding $via$ the Schaeffer correspondence $\Phi$ to the UIPQ, then the labels have a
natural interpretation in terms of distances in the infinite
quadrangulation. In general if  an infinite quadrangulation $q$ is constructed from a labeled tree $\theta=(\tau,\ell)$ in $\mathscr{S}$, every pair $\{u,v\}$ of neighboring vertices in $q$ satisfies $|\ell(u)-\ell(v)| =1$ and thus for every $a,b \in q$ linked by a geodesic $a=a_{0}, a_{1}, \ldots , a_{\op{d}_{\op{gr}}^q(a,b) } =b$ we have the crude bound 
\begin{eqnarray}
\label{trivial}\op{d}_{\op{gr}}^q(a,b) = \sum_{i=1}^{\op{d}_{\op{gr}}^q(a,b) } |\ell(a_{i})-\ell(a_{i-1})|\geq \left|\sum_{i=1}^{\op{d}_{\op{gr}}^q(a,b) } \ell(a_{i})-\ell(a_{i-1}) \right | =  
|\ell(a)-\ell(b)|. 
\end{eqnarray} A better upper bound is given by the  so-called {\em cactus bound}
 \begin{eqnarray} \op{d}^q_{\op{gr}}(a,b)&\geq&  \ell(a)+\ell(b)-2\min_{v\in \llbracket a,b\rrbracket} \ell(v)
 ,  \label{cactus-bound}\end{eqnarray} where we recall that $  \llbracket a,b \rrbracket$ represents the geodesic line in $\tau$ between $a$ and $b$. This bound is   proved in \cite{CLGMcactus} in the context of finite trees and quadrangulations, but remains valid here without change.  The idea goes as follows: let $w$ be of minimal label on $\llbracket a,b\rrbracket$, and assume $w\notin \{a,b\}$ to avoid trivialities.  Removing $w$ breaks the tree $\tau$ into two connected parts, containing respectively $a$ and $b$. Now a path from $a$ to $b$ has to ``pass over'' $w$ using an arc between a corner (in the first component) to its successor (in the other component), and this can only happen by visiting a vertex with label less than $\ell(w)$. Using \eqref{trivial} we deduce that this path at length at least $\ell(a)-\ell(w)+\ell(b)-\ell(w)$, as wanted. 

\subsubsection{From quadrangulations to trees} \label{quadtotrees}
We saw that the Schaeffer mapping $\mathbf{T}_{n}^{(0)}  \times \{0,1\} \longrightarrow \mathbf{Q}_{n}^\bullet$ is in fact a bijection. We now describe the reverse construction. The details can be found in \cite{CS04}. Let $(q,\rho)$ be a finite rooted quadrangulation given with a distinguished vertex $\rho \in V(q)$. We define a labeling $\ell$ of the vertices of the quadrangulation by setting 
 \begin{eqnarray*} \ell(v) &=& \op{d}_{\op{gr}}^q(v,\rho), \ \ v \in V(q). \end{eqnarray*}
Since the map $q$ is bipartite, if $u,v$ are neighbors in $q$ then $|\ell(u)-\ell(v)|=1.$ Thus the faces of $q$ can be decomposed into two subsets: The faces such that the labels of the vertices listed in clockwise order are $(i,i+1,i,i+1)$ for some $i \geq 0$ or those for which these labels are $(i,i+1,i+2,i+1)$ for some $i\geq 0$. We then draw on top of the quadrangulation an edge in each face according to the rules given by the figure below.
\begin{figure}[h] 
\begin{center}
\includegraphics[]{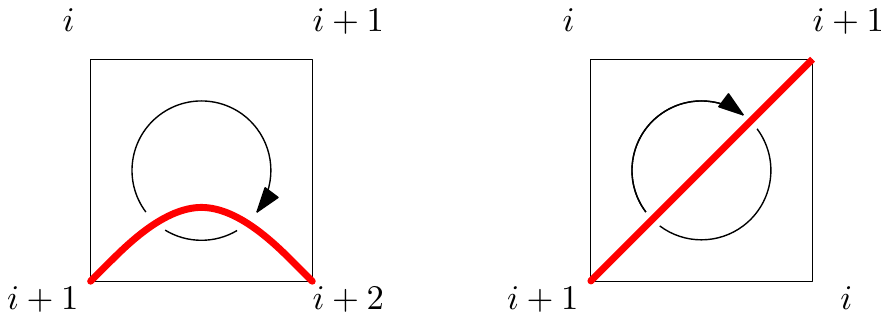}
\caption{\label{rules}Rules for the reverse Schaeffer construction.}
\end{center} 
\end{figure}

\begin{figure}[!h]
\begin{center}
\includegraphics[height=8cm]{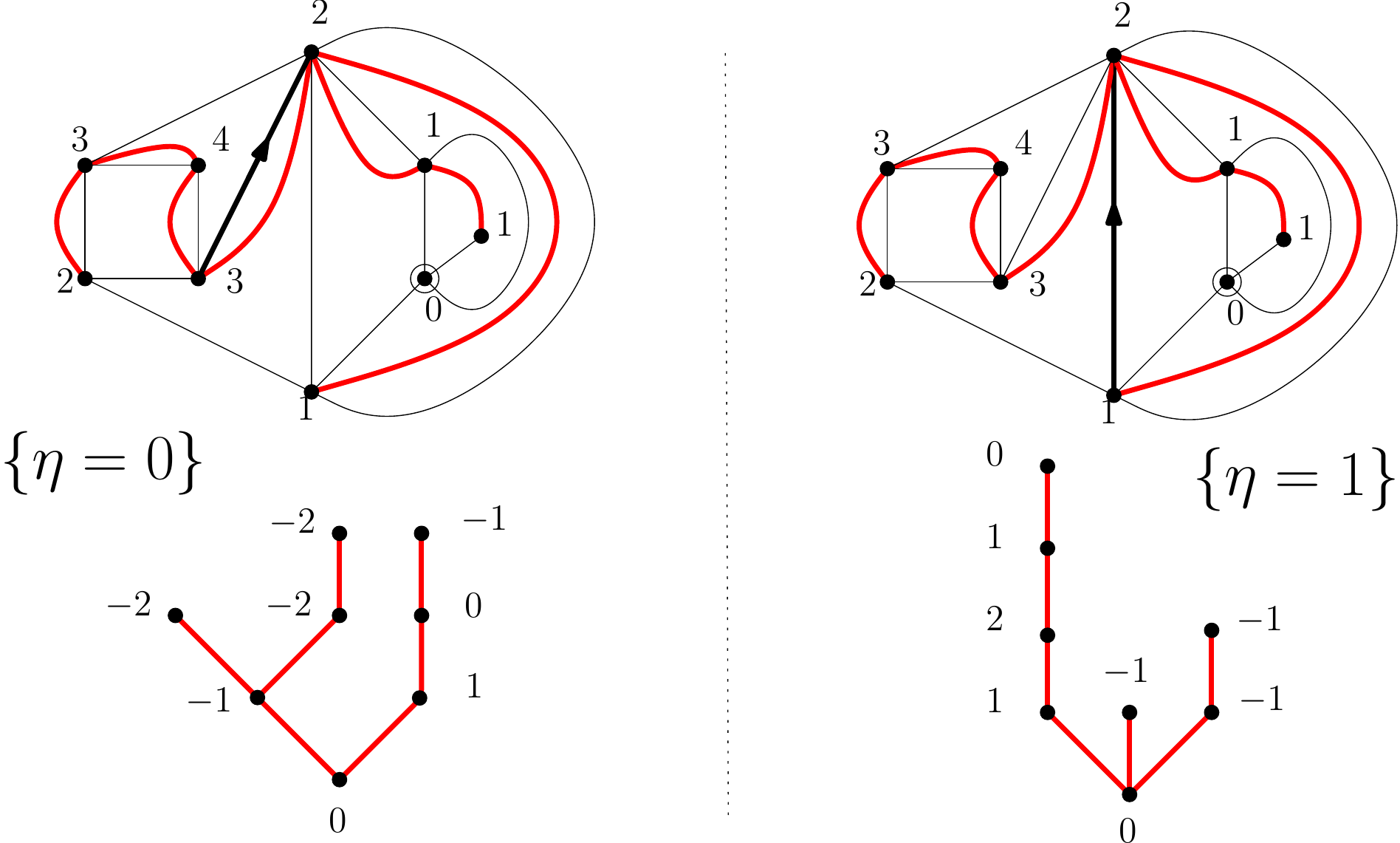}
\caption{ \label{uipq:enracinement} Illustration of the rooting of the plane tree $\tau$}
\end{center}
\end{figure}

The graph $\tau$ formed by the edges added in the faces of $q$ is a spanning tree of $q\backslash\{\rho\}$, see \cite[Proposition 1]{CS04}. This tree comes with a natural embedding in the plane, and we root $\tau$ according to the following rules (see Fig.\ref{uipq:enracinement}):
\begin{itemize}
\item If $\ell(e^*_{-})>\ell(e^*_{+})$ then we root $\tau$ at the corner incident to the edge $e^*$ on $e^*_{-}$, 
\item otherwise we root $\tau$ at the corner incident to the edge $e^*$ on $e^*_{+}$, \end{itemize}
 Finally, we shift the labeling of $\tau$ inherited from the labeling on $V(q)\backslash \{\rho\}$ by the label of the root of $\tau$, 
$$ \tilde{\ell}(u) = \ell(u)-\ell(\varnothing), \ \ u \in \tau.$$
Then we have \cite[Proposition 1]{CS04}  \begin{eqnarray*}
\Phi^{-1}\big((q,\rho)\big) &=& \big((\tau, \tilde{\ell}),\mathbf{1}_{\ell(e^*_{+})>\ell(e^*_{-})}\big). \end{eqnarray*}

\paragraph{Infinite case.} If $q$ is a (possibly infinite) quadrangulation and $\ell: V(q) \longrightarrow \mathbb{Z},$ is a labeling of the vertices of $q$ such that for any neighboring vertices $u,v$ we have $|\ell(u)-\ell(v)|=1$, then a graph can be associated to $(q,\ell)$ by the  device we described above. This graph could contain cycles and is not a tree in the general case.\\
\label{reconstruction}However, suppose that the infinite quadrangulation $q$ is constructed as the image under $\Phi$ of a labeled tree $\theta=(\tau,\ell) \in \mathscr{S}$ and an element of $ \{0,1\}$. Then, with the usual identification of $V(q)$ with $\tau$, the labeling of $V(q)$ inherited from the labeling $\ell$ of $\tau$ satisfies $|\ell(u)-\ell(v)|=1$ for any $u,v \in V(q)$. An easy adaptation of the argument of \cite[Property 6.2]{CD06} then shows that the faces of $q$ are in one-to-one correspondence with the edges of $\tau$ and that the edges constructed on top of each face of $q$ following the rules of Fig. \ref{rules} exactly correspond to the edges of $\tau$. In other words, provided that $q$ is constructed from $\theta=(\tau,\ell)$ then the graph constructed on top of $q$ using the labeling $\ell$ is exactly $\tau$. The rooting of $\tau$ is also recovered from $q$ and $\ell$ by the same procedure as in the finite case.

\subsection{The uniform infinite labeled tree}
\label{sec:uiltuipq}

For every integer $l>0$, we denote by $\rho_l$ the law of the
Galton-Watson tree with geometric offspring distribution with
parameter $1/2$, labeled according to the following rules. The root
has label $l$ and every other vertex has a label chosen uniformly in
$\{m-1,m,m+1\}$ where $m$ is the label of its parent, these choices
being made independently for every vertex. Otherwise said, for every tree
$\theta \in \ltrees{}{l}$, $\rho_l(\theta) = \frac{1}{2}
12^{-|\theta|}$.
\begin{definition} \label{deftinfty}Let $\theta
  =(T_{\infty}, (\ell(u))_{u\in T_{\infty}})$ be a random variable with values in $(\mathbf{T}^{(0)},d_{\arbres})$ whose distribution $\mu$ is described by the following 
properties
\begin{enumerate}
\item $\theta$ belongs to $\mathscr{S}$ almost surely, 
\item the process $(\ell(\mathrm{S}_{T_{\infty}}(n)))_{n \geq 0}$ is a random walk with
  independent uniform steps in $\{-1,0,1\}$,
 \item conditionally given $(\ell(\mathrm{S}_{T_{\infty}}(n)))_{n \geq 0} = (x_n)_{n \geq 0}$, the
   sequence $(L_n(\theta))_{n \geq 0}$ of subtrees of $\theta$ attached to the
   left side of the spine and the sequence $(R_n(\theta))_{n \geq 0}$ of
   subtrees attached to the right side of the spine form two
   independent sequences of independent labeled trees distributed
   according to the measures $(\rho_{x_n})_{n\geq 0}$.
\end{enumerate}
\end{definition}
 In other words, if $\theta=(T_{\infty},\ell)$ is distributed according to $\mu$ then the structure of the tree $T_{\infty}$ is given by an infinite spine and independent critical geometric Galton-Watson trees grafted on the left and right of each vertex of the spine. Conditionally on $T_{\infty}$ the labeling is given by independent variables uniform over $\{-1,0,+1\}$ assigned to each edge of $T_{\infty}$, which represent  the label increments along the different edges, together with the boundary condition $\ell(\varnothing)=0$.
 
 The random infinite tree $T_{\infty}$,  called the critical geometric Galton-Watson tree \emph{conditioned to survive}, was constructed in \cite[Lemma 1.14]{Kes86} as the limit of critical geometric Galton-Watson conditioned to survive up to level $n$, as $n\to \infty$. To make the link between the classical construction of $T_{\infty}$ (see e.g. \cite[Chapter 12]{LP10}) and the one provided by the last definition, note the following equality in distribution 
 \begin{eqnarray*} 1+G+G' &\overset{(d)}{=}& \hat{G},  \end{eqnarray*}
  where  $G,G',\hat{G}$ are independent random variables such that $G,G'$ are geometric of parameter $1/2$ and $\hat{G}$ is a size-biased geometric $1/2$ variable, that is $ \mathbb{P}(\hat{G}=k)=k \P{G=k}=k2^{-(k+1)}$.

The law $\mu$ can also be seen as the law of a uniform infinite element of
$\mathscr{S}$, as formalized by the following statement. 

\begin{thm}[{\cite{Kes86}}]
\label{sec:unif-infin-label}
For every $n \geq 1$, let $\mu_n$ be the uniform probability measure
on $\ltrees{n}{0}$.  Then the sequence $(\mu_n)_{n \in \mathbb{N}}$
converges weakly to $\mu$ in the space of Borel probability measures on $( \ltrees{}{0} ,
d_{\arbres})$.
 \end{thm}
 \proof It is a standard result \cite{LG06} that the distribution of a uniformly chosen planar tree $T_{n}$ with $n$ edges is the same as the distribution of a critical Galton-Watson tree with geometric offspring distribution conditioned on the total progeny to be $n+1$. The convergence in distribution of $T_{n}$ towards $T_{\infty}$ in the sense of $d_{\mathbf{T}}$ then follows from \cite[Lemma 1.14]{Kes86}, see also \cite{LP10}. An analogous result holds for the uniform labeled trees since the labeling is given by independent variables uniform over $\{-1,0,+1\}$ assigned to each edge of the trees.\endproof

  \subsection{The main result}

We are now ready to state our main result. Recall that $\nu$ is the
law of the UIPQ as defined in Theorem \ref{krikun}. Let also
$\mathcal{B}(1/2)$ be the Bernoulli law $(\delta_0+\delta_1)/2$, and
recall the Schaeffer correspondence
$\Phi:\mathscr{S}\times\{0,1\}\to \carte$.  In the following
statement, if $q$ is an element of $\carte_\infty$, and $f:V(q)\to
\mathbb{Z}$ is a function on $V(q)$, we say that $\lim_{z\to\infty}
f(z)=l$ if 
$f$ is equal to $l$ everywhere but
on a finite subset of $V(q)$.

\begin{theorem}\label{sec:unif-infin-label-1}
The probability measure $\nu$ is the image of $\mu\otimes
\mathcal{B}(1/2)$ under the mapping $\Phi$: 
 \begin{eqnarray}\label{eq:3}
\nu &=& \Phi_*\big(\mu\otimes\mathcal{B}(1/2)\big)
 \end{eqnarray}
Moreover, if $(\theta=(T_{\infty},\ell),\eta)$ has distribution
$\mu\otimes\mathcal{B}(1/2)$ and $Q_{\infty}= \Phi(\theta,
\eta)$, then, with the usual identification of the vertices of
$Q_\infty$ with the vertices of $\theta$, one has, almost surely,
 \begin{eqnarray}\label{eq:4}
\ell(u)-\ell(v) &=& \lim_{z \to \infty} \left(
\op{d}_{\op{gr}}^{Q_{\infty}}(u,z) -
\op{d}_{\op{gr}}^{Q_{\infty}}(v,z) \right)\, ,\qquad \forall u,v\in
V(Q_\infty)\, .
 \end{eqnarray}
\end{theorem}

Let us make some comments about this result. The first part of the
statement is easy: Since $\Phi$ is continuous from
$(\mathscr{S}\cup \ltrees{f}{0})\times\{0,1\}$ to $\carte$ and since,
if $\nu_n$ is the uniform law on $\carte_n$, one has
$$\nu_n=\Phi\big(\mu_n\otimes \mathcal{B}(1/2)\big)\, ,$$ and one
obtains (\ref{eq:3}) simply by passing to the limit $n \to \infty$ in this identity using Theorems
\ref{krikun} and \ref{sec:unif-infin-label}. To be completely
accurate, the mapping $\Phi$ in the previous display should be
understood as taking values in $\carte_n$ rather than
$\carte_n^\bullet$, simply by ``forgetting'' the distinguished vertex
arising in the construction of Schaeffer's bijection.

The rest of the statement is more subtle, and says that the labels,
inherited on the vertices of $Q_\infty$ in its construction  from a
labeled tree $(T_{\infty},\ell)$ distributed according to $\mu$, can be recovered as a measurable function
of $Q_\infty$. This is not obvious at first, because a formula such as 
(\ref{eq:2}) is lacking in the infinite setting. It should be replaced
by the asymptotic formula (\ref{eq:4}), which specializes to
\begin{equation}\label{eq:6}
\ell(u)=\lim_{z\to\infty}\big(\op{d}^{Q_{\infty}}_{\mathrm{gr}}(z,u)-
\op{d}^{Q_{\infty}}_{\mathrm{gr}}(z,e^*_\pm)\big)\, ,\, \qquad u\in V(Q_\infty)\,
,
\end{equation}
where
\begin{equation}\label{eq:5}
e^*_\pm=\left\{\begin{array}{cl} e^*_- &\mbox{ if
}\lim_{z\to\infty}(\op{d}^{Q_\infty}_{\mathrm{gr}}(e^*_-,z)-
\op{d}^{Q_\infty}_{\mathrm{gr}}(e^*_+,z))=-1\\ e^*_+ &\mbox{ if
}\lim_{z\to\infty}(\op{d}^{Q_\infty}_{\mathrm{gr}}(e^*_-,z)-
\op{d}^{Q_\infty}_{\mathrm{gr}}(e^*_+,z))=1
\end{array}
\right.\, .
\end{equation}
Of course, the fact that the limits in (\ref{eq:4}) and (\ref{eq:5})
exist is not obvious and is part of the statement. This was first
observed by Krikun in \cite{Kri08}, and will be derived here by 
different methods. Note that the vertex $e^*_\pm$ corresponds to the root vertex $\varnothing$ of $T_{\infty}$ in the natural identification of vertices of $Q_{\infty}$ with vertices of $T_{\infty}$.

In particular, the fact that the labels are measurable with respect to
$Q_\infty$ entails that $(\theta,\eta)$ can be recovered as a measurable function of $Q_{\infty}$. Indeed, by the discussion at the end of Section \ref{quadtotrees}, the tree $T_{\infty}$ can be reconstructed from $Q_{\infty}$ and the labeling $\ell$. The Bernoulli variable $\eta$ is also recovered by \eqref{eq:5}. This settle the three conjectures proposed by Krikun in \cite{Kri08}.
\medskip


The proof of (\ref{eq:4}) depends on certain properties of geodesics in the UIPQ that we derive in the next section. Before this, we give another view on our result in terms of asymptotic geometry of the UIPQ.

\subsubsection{Gromov compactification of the UIPQ}\label{sec:grom-comp-uipq} 

Let $(X,d)$ be a locally compact metric space. The set $C(X)$ of real-valued continuous functions on $X$ is endowed with the topology of uniform convergence on every compact set of $X$. One defines an equivalence relation on $C(X)$ by declaring two functions equal if they differ by an additive constant and  the associated quotient space endowed with the quotient topology is denoted by $C(X)/\mathbb{R}$. Following \cite{Gro81}, one can embed the original space $X$  in $C(X)/\mathbb{R}$ using the injection
$$ i : \begin{array}{ccccc} 
 X &\longrightarrow&{C}(X)& \longrightarrow& {C}(X)/\mathbb{R}\\ x&\longmapsto & d_{x}= d(x,.) & \longmapsto & \overline{d_{x}}\, .
\end{array}
$$
The {\em Gromov compactification} of $X$ is then the closure of $i(X)$ in $C(X)/\mathbb{R}$. The Gromov boundary $\partial X$ of $X$ is composed of the points in the closure of $i(X)$ in $C(X)/ \mathbb{R}$ which are not already in $i(X)$. The points in $\partial X$ are called \emph{horofunctions}, see \cite{Gro81}.

Applying this discussion to the case where $(X,d)=(V(Q_\infty),\op{d}_{\op{gr}}^{Q_\infty})$, we can immediately rephrase the last part of Theorem \ref{sec:unif-infin-label-1} as follows. 

\begin{corollary} Almost surely, the Gromov boundary $\partial Q_{\infty}$ of the UIPQ consists of only one point which is $\overline{\ell}$, the equivalence class of $\ell$ up to additive constants.
\end{corollary}

\section{Geodesics in the UIPQ}\label{sec:geodesics-uipq}

\paragraph{Geodesics.} If $G=(V,E)$ is a graph, a {\em chain} or \emph{path} in $G$ is a (finite or
infinite) sequence of vertices $\gamma = (\gamma(0),\gamma(1),
\ldots)$ such that for every $i\geq 0$, the vertices $\gamma(i)$ and
$\gamma(i+1)$ are linked by an edge of the graph. Such a chain is called a
\emph{geodesic} if for every $i,j \geq 0$, the graph distance $\op{d}_{\op{gr}}^G$ between
$\gamma(i)$ and $\gamma(j)$ is equal to $|j-i|$.  A {\em geodesic ray} emanating from $x$ is an infinite geodesic starting at $x\in V$. 
\bigskip	

We will establish two properties of the geodesics in the UIPQ: A confluence property towards the root  (Section \ref{sec:confluent-geodesics}) and a confluence property towards infinity (Section \ref{sec:infinite-geodesics}). These two properties are reminiscent of the work of Le Gall on geodesics in the Brownian Map \cite{LG09}. Put  together they yield the last part \eqref{eq:4} of Theorem \ref{sec:unif-infin-label-1}.

\subsection{Confluent geodesics to the root}\label{sec:confluent-geodesics}

Let $Q_{\infty}$ be distributed according to $\nu$ (see Theorem \ref{krikun}) and $x$ be a vertex in $Q_\infty$. For every  $R\geq 0$, we want to show that (with probability $1$)
it is possible to find $R'\geq R$ and a family of geodesics
$\gamma^z_R,z\notin B_{\carte,R'}(Q_\infty)$ linking $x$ to $z$
respectively, such that for every $z,z'\notin
B_{\carte,R'}(Q_\infty)$,
$$\gamma^z_R(i)=\gamma^{z'}_R(i)\, ,\qquad \mbox{ for every }i\in
\{0,1,\ldots,R\}\, .$$ In other words, all of these geodesics start
with a common initial segment, independently of the target vertex $z$.

To this end, we need the construction by Chassaing-Durhuus \cite{CD06} of the UIPQ, which we briefly recall. Let $l\geq 1$ and set $\wltrees{}{l}$ be the subset of $\ltrees{}{l}$ consisting of all trees $\theta=(\tau,\ell)$ such that $\ell(v)\geq 1$ for every $v\in \tau$. Elements of $\wltrees{}{l}$ are called $l$-well-labeled trees, and just well-labeled trees if $l=1$.  We let $\wltrees{n}{l}$ (resp.\ $\wltrees{\infty}{l}$) be the set of all $l$-well-labeled trees with $n$ edges (resp.\ of infinite $l$-well-labeled trees).

Let $\overline{\mu}_n$ be the uniform distribution on
$\overline{\mathbf{T}}\ \hspace{-0.9mm}_{n}^{(1)}$. Let also $\overline{\mathscr{S}}$ be the set of all trees
$\theta=(\tau,\ell)\in \overline{\mathbf{T}}\ \hspace{-0.9mm}_{\infty}^{(1)}$ such that 
\begin{itemize}
\item
the tree $\tau$ has a unique spine, and
\item
for every $R\in
\mathbb{N}$, the set $\{v\in \tau:\ell(v)\leq R\}$ is finite. 
\end{itemize}

\begin{proposition}[\cite{CD06}]
\label{sec:confluent-geodesics-1}
The sequence $(\overline{\mu}_n)_{n \geq 1}$ converges weakly to a limiting
probability law $\overline{\mu}$, in the space of Borel probability
measures on $(\wltrees{}{1},d_{\arbres})$. Moreover, we have 
$\overline{\mu}(\overline{\mathscr{S}})=1$.
\end{proposition}

The exact description of $\overline{\mu}$ is not important for our
concerns, and can be found in \cite{CD06}.  The Schaeffer correspondence $\overline{\Phi}$ can be defined on $\overline{\mathscr{S}}$. Let us describe quickly this correspondence. Details can be found in \cite{CD06}, see also \cite{LGM10,Men08}.

Let $\theta=(\tau,\ell)$ be an
element of $\overline{\mathscr{S}}$. We start by embedding $\tau$ in the plane
in such a way that there are no accumulation points (which is possible
since $\tau$ is locally finite). We add an extra vertex $\partial$ in the plane, not belonging to the embedding of $\tau$. Then, we let
$(c_i^{(L)},i\geq 0)$ and $(c_i^{(R)},i\geq 0)$ be the sequence of
corners visited in contour order on the left and right sides, starting
with the root corner of $\tau$. We let, for $i\in \mathbb{Z}$,
$$c_{i}=\left\{\begin{array}{cc}
c_i^{(L)} & \mbox{ if } i\geq 0\\
c_{-i}^{(R)} & \mbox{ if } i< 0
\end{array}\right. .$$ 
We now define the notion of successor. If the label of $\mathcal{V}(c_i)$
is $1$, then the successor of the corner $c_i$ is
$\partial$. Otherwise, the successor of $c_i$ is the first corner
$c_j$ in the infinite list $\{c_{i+1},c_{i+2},\ldots\} \cup \{ \ldots , c_{i-2},c_{i-1}\}$
such that $\ell(c_j)=\ell(c_i)-1$. The successor of any corner $c_i$
with $\ell(c_i)\geq 2$ exists because of the labeling constraints, and
the definition of $\overline{\mathscr{S}}$.

The end of the construction is similar to Section \ref{sec:schaeffer}: We draw an edge between
each corner and its successor and then remove all the edges of the
embedding of $\tau$. The new edges can be drawn in such a way that the resulting embedded graph is proper and represents an infinite quadrangulation of the plane. We denote this quadrangulation by  $\overline{\Phi}(\theta)$ and
root it at the arc from $\partial$ to $c_0$. Note that in this
construction, we do not need to introduce an extra parameter
$\eta$ to determine the orientation of the root. Moreover the non-negative labels $\ell$ have the following interpretation in terms of distances in $\overline{\Phi}(\theta)$. For every $u \in \tau$,
\begin{eqnarray} \ell(u) &=& \op{d}_{\op{gr}}^{\overline{\Phi}(\theta)}(\partial,u), \label{dist:positif}\end{eqnarray} with the identification of the vertices of $\overline{\Phi}(\theta)$ with $\tau \cup \{\partial\} $.

\begin{proposition}[\cite{CD06},\cite{Men08}]\label{sec:confluent-geodesics-2}
It holds that 
 \begin{eqnarray*}\nu&=&\overline{\Phi}_*\overline{\mu}\, ,  \end{eqnarray*} that is, the UIPQ  follows the distribution of
$\overline{\Phi}(\theta)$, where $\theta$ is random with distribution
$\overline{\mu}$.
\end{proposition}

Notice that the mapping $\overline{\Phi}:\overline{\mathscr{S}}\to \carte$ is injective. Its inverse function $\overline{\Phi} {}^{-1}:\overline{\Phi}(\overline{\mathscr{S}})\to \wltrees{}{1}$ is described in a similar manner as in Section \ref{quadtotrees}: Given the quadrangulation $q= \overline{\Phi}(\tau,\ell)$, we recover the labeling $\ell$ over $V(q)\backslash \{\partial\}$ by \eqref{dist:positif} and $\ell(\partial)=0$. Note that $\partial$ is always the origin of the root edge of $q$. We then apply the same device as for $\Phi^{-1}$, that is, separating the faces of $q$ into two kinds and adding an edge on top of them according to Fig. \ref{rules}. The resulting graph is $\tau$ and is rooted at the corner incident to the root edge of $q$. One can check that the mapping $\overline{\Phi} {}^{-1}$ is continuous, i.e.\ that for every $h>0$, the neighborhood $B_{\mathbf{T},h}(\tau,\ell)$ is determined by $B_{\mathbf{Q},r}(q)$ as soon as $r$ is large enough.  Thus if $Q_{\infty}$ is distributed according to $\mu$, one can define a labeled tree $(\tau,\ell)$ distributed according to $\overline{\mu}$ as a measurable function of $Q_{\infty}$ such that $Q_{\infty} = \overline{\Phi}(\tau,\ell)$.
\medskip

\begin{figure}[!h]
 \begin{center}
 \includegraphics[width=16cm]{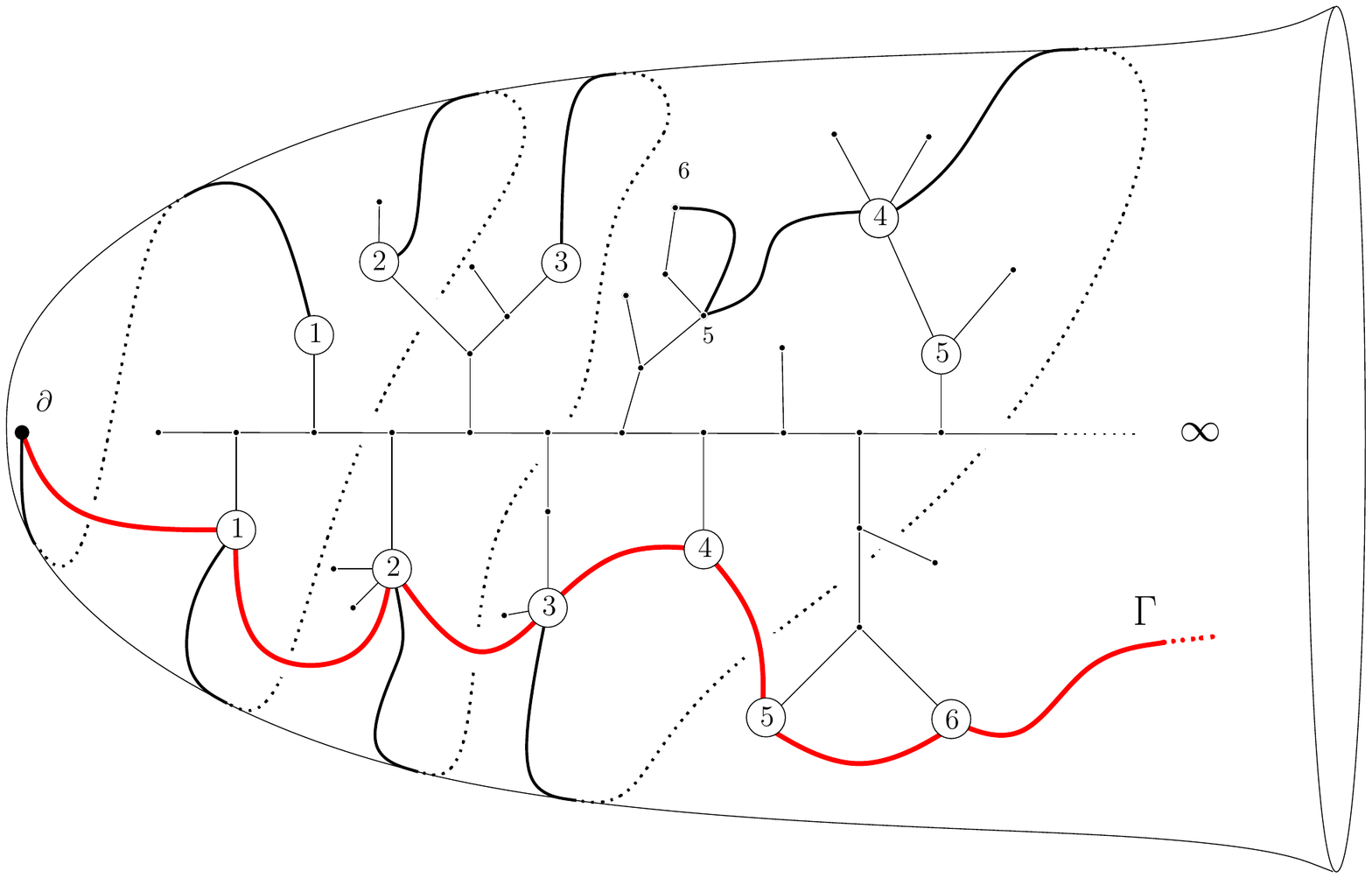}
 \caption{ \label{schpositif}Illustration of the proof of Lemma \ref{sec:confluent-geodesics-3}. The tree is represented in solid lines. Every vertex marked by a circled integer corresponds to the last occurrence of this integer along either the left or the right  side of the tree. The distinguished geodesic $\Gamma$ is represented by a thick line.}
 \end{center}
 \end{figure}

From this construction, it is possible to specify a particular
infinite geodesic (or {\em geodesic ray}) starting from $e^*_-$.   Namely, if $(c_{i})_{i\in \mathbb{Z}}$ is the contour sequence of $\tau$, for
every $i\geq 1$, let
 \begin{eqnarray*} d(i) &=&\min\{j\leq 0:\ell(c_{j})=i\}\, , \end{eqnarray*} which is finite by definition of $\overline{ \mathscr{S}}$. Then there is an arc
between $c_{d(i+1)}$ and $c_{d(i)}$ for every $i\geq 1$, as well as an
arc from $c_{d(1)}$ to $\partial$, and the path $(\partial,
\mathcal{V}(c_{d(1)}),\mathcal{V}(c_{d(2)}),\ldots)$ is a geodesic ray. We
call it the {\em distinguished geodesic ray} of $Q_{\infty}$, and denote it by
$\Gamma$, see Fig.\,\ref{schpositif}. 
 
\begin{lemma}\label{sec:confluent-geodesics-3}
For every $R\geq 0$, there exists $R'\geq R$ such that every $z\in
V(Q_\infty)\setminus B_{\carte,R'}(Q_\infty)$ can be joined to
$\partial$ by a geodesic chain $\gamma$ such that
$\gamma(i)=\Gamma(i)$ for every $i\in \{0,1,2,\ldots,R\}$. 
\end{lemma}

\proof Let $Q_\infty$ be distributed according to $\nu$ and set $(\tau,\ell) = \overline{\Phi}{}^{-1}(Q_\infty)$. Finally define $\Gamma$ as above. Define
$$R'=\max_{d(R)\leq i\leq g(R)} \ell(c_i)\, ,$$
where $d(R)$ is defined above, and 
$$g(i)=\max\{j\geq 0:\ell(c_{j})=i\}\, .$$ Let $z$ be a vertex of
$Q_\infty$, not in $B_{\carte,R'}(Q_\infty)$, and let $c_j$ be any
corner incident to $z$. Then $j$ cannot be in $[d(R),g(R)]$ since by
definition $\ell(c_j)=\op{d}^{Q_\infty}_{\mathrm{gr}}(\partial,z)>R'\geq
\ell(c_i)$ for any $i\in [d(R),g(R)]$. Now, let $\gamma$ be the
geodesic defined as the path starting at $c_j$, and following the arcs
from $c_j$ to its successor corner, then from this corner to its
successor, and so on until it reaches $\partial$. These
geodesics have the desired property, see Fig.\,\ref{schpositif}. Note that if $j>0$, that is, if
$c_j$ lies on the left side of $\tau$, then necessarily all vertices in
the geodesic $\gamma$ with label less than or equal to $R$ have to lie
on the right-hand side of $\tau$. See Fig.\,\ref{schpositif}. \endproof

\subsection{Coalescence of proper geodesics rays to infinity}
\label{sec:infinite-geodesics}

With the notation of Theorem \ref{sec:unif-infin-label-1}, let $(\theta=(T_{\infty},\ell),\eta)$ be distributed according to $\mu\otimes \mathcal{B}(1/2)$, and let $Q_\infty$ be the image of $(\theta,\eta)$ by the Schaeffer correspondence $\Phi$. The construction of $Q_\infty$ from a tree $\theta \in  \mathscr{S}$ allows to specify another class of geodesic rays in $Q_\infty$, which are defined as follows. These geodesic rays are emanating from the root vertex $\varnothing$ of $\theta$, which can be either $e^*_-$ or $e^*_+$, depending on the value of $\eta$. Consider any infinite path $(u_0,u_1,u_2,\ldots)$ in $q=\Phi(\theta,\eta)$ starting from $\varnothing=u_0$, and such that $\ell(u_i)=-i$ for every $i$. Then necessarily, such a chain is a geodesic ray emanating from $\varnothing$, because from (\ref{trivial}) we have $\op{d}^{q}_{\mathrm{gr}}(u_i,u_j)\geq |i-j|$ for every $i,j\geq 0$, and the other inequality is obviously true.

We call such a geodesic a {\em proper} geodesic ray emanating from $\varnothing$. We will see in Corollary \ref{sec:coal-geod-rays} below that all geodesic rays emanating from $\varnothing$ are in fact proper. 

It should be kept in mind that in the definition of a geodesic ray $(\gamma(i),i\geq 0)$, we can further specify which edge of $Q_\infty$ is used to pass from $\gamma(i)$ to $\gamma(i+1)$. If $\gamma$ is a proper geodesic, then $\ell(\gamma(i+1))=\ell(\gamma(i))-1$ so this edge has to be an arc drawn in the Schaeffer correspondence between a corner $c^\gamma(i)$ of $T_\infty$ and its successor $\mathcal{S}(c^\gamma(i))$. We will use this several time in the sequel. 

The main result of this section shows the existence of cut-points visited by  every infinite proper geodesic.  \begin{theorem}\label{cut-points}Let $(\theta=(T_{\infty},\ell), \eta)$ be distributed according to $\mu \otimes \mathcal{B}(1/2),$ and let $Q_{\infty}$ be the image of $(\theta,\eta)$ by the Schaeffer correspondence $\Phi$. Almost surely, there exists an infinite
  sequence of distinct vertices $p_1,p_2,\ldots \in V(Q_\infty)$ such that every proper geodesic ray emanating from $\varnothing$ passes through $p_1,p_2,\ldots$.  \end{theorem}
  
  To prove this theorem we will introduce two specific proper geodesic rays that are in a sense extremal: the minimal and maximal geodesics. We then prove that they meet infinitely often (Lemma \ref{sec:proof-theorem-refcut}) and that these meetings points eventually are common to all proper geodesics emanating from $\varnothing$. 

\subsubsection{The maximal and minimal geodesics}

 Recall that if $\theta=(\tau,\ell) $ is a labeled tree in $\mathscr{S}$ with contour sequence $(c_{i})_{i \in \mathbb{Z}}$, for every $j \in \mathbb{Z}$ the successor $\mathcal{S}(c_{j})$ of $c_{j}$ is the first corner among $c_{j+1}, c_{j+2}, \ldots$ with label $\ell(c_{j})-1$. 
\begin{definition}[maximal geodesic] \label{def:maxgeo} Let $\theta=(\tau,\ell) \in \mathscr{S}$. For every corner $c$ of $\theta$, the \emph{maximal geodesic} $\gamma_{\max}^c$ emanating from $c$ in $\theta$ is given by the chain of vertices attached to the iterated successors of $c$, 
$$ \gamma_{\max}^c(i) := \mathcal{V}\big(\mathcal{S}^{(i)}(c)\big), \ \ i \geq 0,$$ where $\mathcal{S}^{(i)}$ is the $i$-fold composition of the successor mapping. 
\end{definition}
Using \eqref{trivial} again, we deduce that the maximal geodesics are indeed geodesic chains in the quadrangulation associated to $\theta$. When $c=c_{0}$ is the root corner of $\tau$ we drop $c_{0}$ in the notation $\gamma_{\max}$ and call it the maximal geodesic. The maximal geodesic is a proper geodesic, and in the above notation, $c^{\gamma_{\max}}(i)=\mathcal{S}^{(i)}(c_0)$. 

Next, consider only the left part of an infinite labeled tree $(\tau, \ell)$, which corresponds to the corners $(c_{0},c_{1},c_{2}, \ldots)$.  We define the \emph{minimal geodesic} $\gamma_{\min}$ inductively. First, let $\gamma_{\min}(0)=\varnothing$. Suppose that the first $n$ steps $(\gamma_{\min}(0), \ldots , \gamma_{\min}(n))$ of $ \gamma_{\min}$ have been constructed. We then set $c^{\gamma_{\min}}(n)$ to be the last corner among $c_{0},c_{1}, \ldots$ that is incident to the vertex $\gamma_{\min}(n)$, which implies 
 \begin{eqnarray*} \gamma_{\min}(n+1)  &=& \mathcal{V}\big( \mathcal{S}(c^{\gamma_{\min}}(n))\big).  \end{eqnarray*} One can check by induction that   $\ell(\gamma_{\min}(i))=-i$, thus $\gamma_{\min}$ is a   proper geodesic ray emanating from $\varnothing$ in $q$.
We restrict the definition of the minimal geodesic to the left part of the tree in order to prescribe the behavior of the path when it hits the spine of the tree. Roughly speaking, the minimal geodesic can hit the spine of $\tau$, but it cannot cross it. 

\begin{figure}[h]
\begin{center}
\includegraphics[width=15cm]{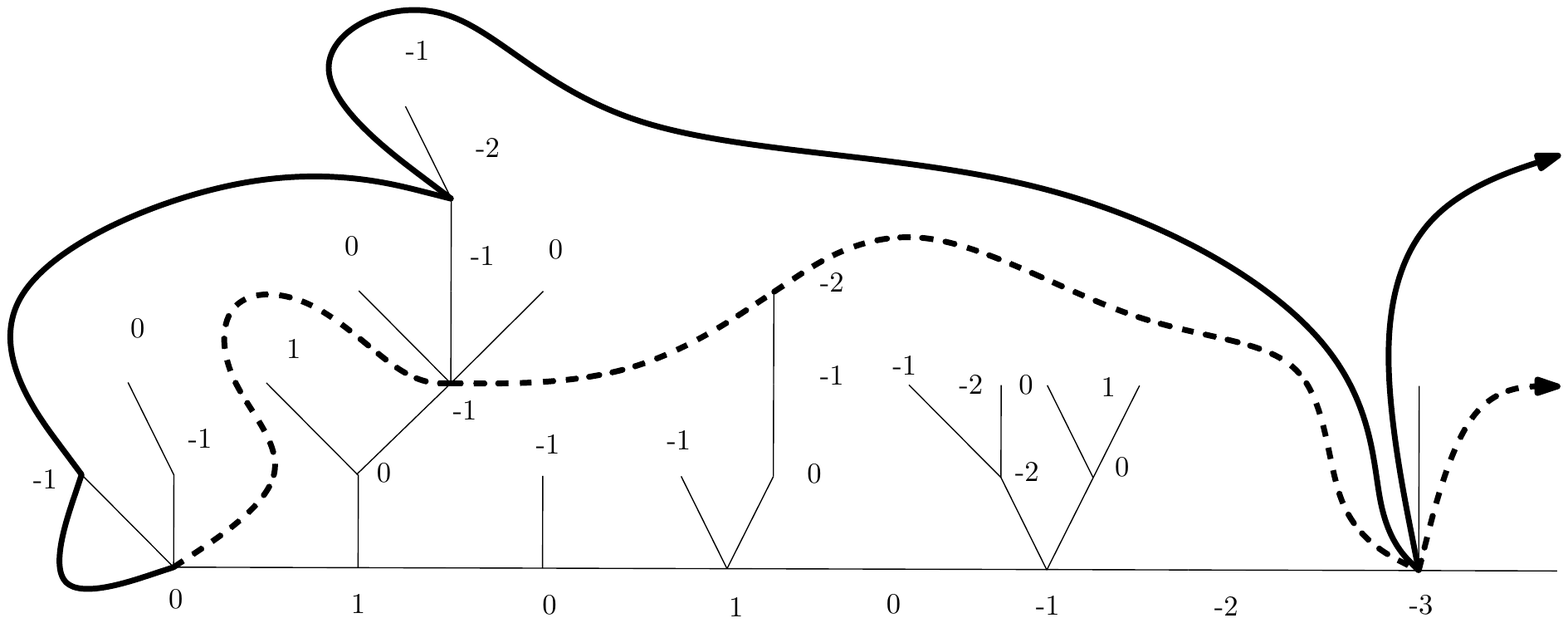}
\caption{The maximal (in solid line) and minimal (in dotted line) geodesics starting from the root corner of the tree $\theta$.}
\end{center}
\label{fig:twogeo}
\end{figure}

The next geometric lemma roughly says that any proper geodesic $\gamma$ is stuck in between $\gamma_{\min}$ and $\gamma_{\max}$ except when $\gamma_{\min}$ hits the spine of $T_{\infty}$ in which case $\gamma$ can visit the right part of the tree $T_{\infty}$.

\begin{lemma}\label{sec:maxim-minim-geod-1} 
Let $\gamma$ be a proper geodesic. 
\begin{enumerate} \item[{\rm (i)}] Suppose that $c$ is a corner incident to $\gamma(i)$ that lies on the left-hand side of $\tau$. Then $c^{\gamma_{\max}}(i)\leq c\leq c^{\gamma_{\min}}(i)$.  
\item[{\rm (ii)}] For every $i\geq 0$, the vertices 
$\gamma(i+1),\gamma(i+2),\ldots$ do not belong to $\llbracket \varnothing, \gamma(i)\rrbracket$.  
\item[{\rm (iii)}] For every $i\geq 0$, if $\gamma(i)$ is incident to the left-hand side of $\tau$, then there exists a (unique) $j\leq i$ such that $\gamma_{\min}(j)$ is an ancestor of $\gamma(i)$ in its subtree to the left of the spine. This means that $\llbracket \varnothing,\gamma(i)\rrbracket$ contains $\gamma_{\min}(j)$, but $\llbracket\gamma_{\min}(j),\gamma(i)\rrbracket\setminus \{\gamma_{\min}(j)\}$ does not interset the spine of $\tau$.  \end{enumerate} \end{lemma}

\proof We first prove that for every $v\in \llbracket\varnothing,\gamma_{\min}(i)\rrbracket\setminus \{\gamma_{\min}(i)\}$ in $\tau$, it holds that $\ell(v)>-i$. By definition of the successor, we have $\ell(c)\geq -i$ for every corner $c$ with $c^{\gamma_{\min}}(i)\leq c <\mathcal{S}(c^{\gamma_{\min}}(i))$. But every vertex in $\llbracket \gamma_{\min}(i),\gamma_{\min}(i+1)\rrbracket\setminus \{\gamma_{\min}(i+1)\}$ is incident to a corner as above. We then argue that $\llbracket\varnothing,\gamma_{\min}(i)\rrbracket\setminus \{\gamma_{\min}(i)\}$ is contained in the union of $\llbracket \gamma_{\min}(j),\gamma_{\min}(j+1)\rrbracket\setminus \{\gamma_{\min}(j+1)\}$ for $0\leq j\leq i-1$. 

Next, let $\gamma$ be a proper geodesic and $c$ be a corner incident to $\gamma(i)$.  Since $c^{\gamma_{\max}}(i)$ is the first corner on the left-hand side of $\tau$ with label $-i$, and since $\ell(c)=\ell(\gamma(i))=-i$, if $c\geq c_0$ we must have $c^{\gamma_{\max}}(i)\leq c$. The inequality $c\leq c^{\gamma_{\min}}(i)$ is true even if $c<c_0$, and is proved by induction. For $i=0$ it is obvious that any corner incident to $\gamma(0)=\varnothing$ is less than $c^{\gamma_{\min}}(0)$. Suppose that $c\leq c^{\gamma_{\min}}(i)$ for every $c$ incident to $\gamma(i)$. Then this holds in particular for $c=c^{\gamma}(i)$, and we deduce $\mathcal{S}(c^\gamma(i))\leq \mathcal{S}(c^{\gamma_{\min}}(i))$ because $\mathcal{S}$ is non-decreasing when restricted to the (ordered) set of corners of $\tau$ with label $-i$. Since $c^{\gamma_{\min}}(i+1)$ is the largest corner incident to $\mathcal{V}(\mathcal{S}(c^{\gamma_{\min}})(i))$, we obtain that any corner $c'$ incident to $\gamma(i+1)$ satisfies $c'\leq c^{\gamma_{\min}}(i+1)$, unless $\gamma(i+1)$ is a strict ancestor of $\gamma_{\min}(i+1)$ in the subtree to the left of $\tau$ that contains $\gamma_{\min}(i+1)$. But by (i), every such ancestor of $\gamma_{\min}(i+1)$ has label at least $-i$, while $\gamma(i+1)$ has label $-i-1$ because $\gamma$ is a proper geodesic, so the latter obstruction does not occur.  This proves (i). 

For (ii), since $\ell(\gamma(i))=-i$, if $\gamma(i+j)\in \llbracket\varnothing,\gamma(i)\rrbracket$ for some $j>0$, then the cactus bound  \eqref{cactus-bound} would imply $\op{d}^q_{\op{gr}}(\varnothing,\gamma(i))\geq 0-i-2(-i-j)=i+2j$, which is impossible since $\gamma$ is a geodesic and $\varnothing=\gamma(0)$. 

We finally prove (iii). It follows from the construction of $\gamma_{\min}$ that the wanted property is hereditary in $i$ as long as $c^{\gamma}(i)\geq c_0$. In particular, the property holds for $\gamma=\gamma_{\max}$. To show that it is true for any proper geodesic $\gamma$, the only problem is when $\gamma$ visits the right-hand side of $\tau$, which can happen only after a time $i$ where $\gamma$ hits the spine. So suppose that $\gamma(i)$ belongs to the spine and that $c=c^{\gamma}(i)\leq c_0$. Since $\gamma$ moves by taking successors, the sequence of corners $c^{\gamma}(j),j\geq i$ can only increase as long as they stay less than $c_0$, and after finitely many steps this sequence must leave the right-hand side of $\tau$. In the meantime, it cannot visit the spine, because of (ii). Therefore, at the first time $j$ after $i$ that $\gamma(i)$ leaves the right-hand side of $\tau$, it must make a step from $c^{\gamma}(j-1)\leq 0$ to $\mathcal{S}(c^{\gamma}(j-1))\geq c_0$. Necessarily, this implies that this successor is also the first corner with label $-j$ on the left-hand side of $\tau$, i.e.\ $\gamma(j)=\gamma_{\max}(j)$. Being a point of the maximal geodesic, we already noticed that property (iii) holds at this stage, from which we are back to a hereditary situation until the next time $k$ where $c^\gamma(k)\leq c_0$, which might happen only at a point where $\gamma_{\min}(k)$ is on the spine. We conclude that the property of (iii) is true at every $i$ where $\gamma(i)$ is  incident to the left-hand side, as wanted. 
\endproof

Suppose now that for $i \geq 1$ we have $\gamma_{\min}(i) = \gamma_{\max}(i)$.  We claim that in fact $$ \mathcal{S}(c^{\gamma_{\max}}(i-1) ) = \mathcal{S}(c^{\gamma_{\min}}(i-1) ).$$ Indeed by property (i) we have $c^{\gamma_{\max}}(i-1) \leq c^{\gamma_{\min}}(i-1)$, and since these two corners have same label $-i+1$ this implies $ \mathcal{S}(c^{\gamma_{\max}}(i-1) ) \leq \mathcal{S}(c^{\gamma_{\min}}(i-1) )$. The inequality can be strict only if $\gamma_{\min}(i) = \gamma_{\max}(i)$ is an ancestor of $\gamma_{\min}(i-1)$ which is prohibited by property (ii). So if $\gamma$ is a geodesic such that $ c^\gamma(i-1)$ lies on the left of $T_{\infty}$ then $c^{\gamma_{\max}}(i-1) \leq c^\gamma(i-1) \leq c^{\gamma_{\min}}(i-1)$ by property (i), which forces $ \mathcal{S}(c^{\gamma_{\max}}(i-1) )= \mathcal{S}(c^{\gamma}(i-1) ) = \mathcal{S}(c^{\gamma_{\min}}(i-1) )$. In particular $\gamma_{\max}(i) = \gamma_{\min}(i) = \gamma(i)$.

We are going to show that $\gamma_{\max}(i)=\gamma_{\min}(i)$ happens infinitely often, and that any proper geodesic visits only vertices to the left of $\tau$ eventually, which will entail Theorem \ref{cut-points}.

\subsubsection{Inbetween $\gamma_{\min}$ and $\gamma_{\max}$}

For $i \geq 0$, denote by $A_{i}$ the labeled tree consisting of $\gamma_{\min}(i)$ and its descendants in the left side of the tree $\tau$ (recall that $\gamma_{\min}(i)$ might be on the spine). Note that this tree may consist only of the single labeled vertex $\gamma_{\min}(i)$. It is clear by construction that $A_i$ is an element of $\ltrees{f}{-i}$. Consider the subtree $\mathsf{C}$ of the left-hand side of $\tau$ obtained by chopping off the trees $A_i,i\geq 0$, i.e.\ $v\in \mathsf{C}$ if and only if $v$ is incident to the left-hand side of $\tau$ and for every $i\geq 0$, $v$ is not a strict descendent of $\gamma_{\min}(i)$ in the subtree to the left of $\tau$ that contains the latter. The tree $\mathsf{C}$ inherits a labeling from $\ell$, so it should be seen as an element of $\mathscr{S}$.

\begin{lemma} \label{indpt} Under the law $\mu$, the sequence  $A_{0}, A_{1},    \ldots $ has law $\displaystyle\bigotimes_{i=0}^\infty \rho_{-i}$ and is independent of $\mathsf{C}$. 
\end{lemma}

\begin{figure}[!h]
 \begin{center}
 \includegraphics[width=8cm]{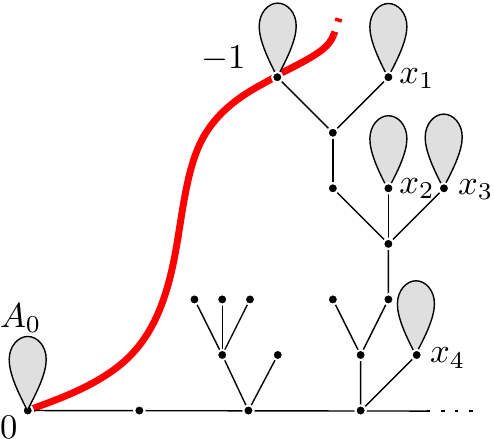}
 \caption{ \label{exploration}Exploration of $T_{\infty}$ until the first $-1$. The gray trees are unexplored.}
 \end{center}
 \end{figure}

\proof[Proof of Lemma \ref{indpt}] First of all, remark that $A_{0}$ is  the tree $L_{0}$ grafted to the left-hand side of the first vertex $ \mathrm{S}(0)$ of the spine of $T_{\infty}$. This tree is shortcut by the minimal geodesic $\gamma_{\min}$ who seeks the successor of $ c^{\gamma_{\min}}(0)$ which is the last corner associated to $\varnothing$. We thus discover the remaining tree $T_{\infty} \backslash A_{0}$ step by step in a depth-first search manner by revealing the children of the vertex $1$ then the children of $11$ (if $11 \in T_{\infty}$) and so on and so forth in lexicographical order. During this exploration, one can obviously discover the labeling in the same time by sampling at each newly discovered edge a uniform variable $\in \{ -1,0,+1\}$ independently of the past carrying the variation of the label along this edge. We then stop when we  discover the first vertex (in lexicographical order) with label $-1$, see Fig.\,\ref{exploration}. 

 This vertex is obviously $\gamma_{\min}(1)$. So far, we thus have explored the part of $T_{\infty}$ which is composed of the vertices on the left of the segment $\llbracket \varnothing, \gamma_{\min}(1)\rrbracket$ together with all the children of the vertices lying on the ancestral path that link $\gamma_{\min}(1)$ to the spine of $T_\infty$. These vertices are denoted by $\gamma_{\min}(1) = x_{0}, x_{1}, x_{2}, \ldots$ in lexicographical order. Note that some of these vertex can have a label equal to $-1$ but none of them has a label strictly less than $-1$. By standard properties of Galton-Watson tree, the subtrees above $x_{0}, x_{1}, x_{2}, \ldots$ are independent of $A_{0}$ and of the part of $T_{\infty}$ explored so far, and form a sequence of independent labeled Galton-Watson trees with laws $\rho_{\ell(x_{i})}$. The tree above $x_{0}= \gamma_{\min}$ is the tree $A_{1}$ that is now shortcut by $\gamma_{\min}$ who seeks the successor of $c^{\gamma_{\min}}(1)$ (which is the last corner incident to $\gamma_{\min}(1)$). This vertex must lie in the unexplored part of $T_{\infty} \backslash (A_{0} \cup A_{1})$. We then continue the exploration of $T_{\infty} \backslash (A_{0} \cup A_{1})$ starting with the tree above ${x_{1}}$ in search of the first vertex (for the lexicographical order) with label $-2$. The process can be carried out iteratively and yields that $A_{0}, A_{1}, A_{2}, \ldots$ has law $\bigotimes_{i \geq 0} \rho_{-i}$ and is independent of $T_{\infty} \backslash \bigcup_{i \geq 0} A_{i}$. 
 \begin{figure}[!h]
  \begin{center}
  \includegraphics[width=8cm]{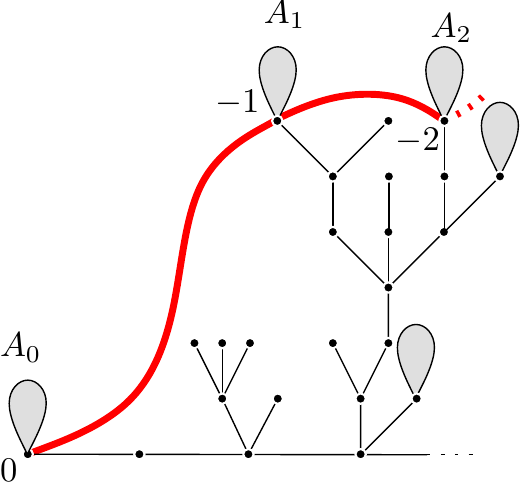}
  \caption{The second step of the exploration.}
  \end{center}
  \end{figure}\endproof

A key step towards Theorem \ref{cut-points} concerns the intersection points between $\gamma_{\max}$ and $\gamma_{\min}$. Let
 \begin{eqnarray*} \mathcal{R}&=&\{i\geq 0:\gamma_{\max}(i)=\gamma_{\min}(i)\},  \end{eqnarray*}
and note that the set $\{\gamma_{\max}(i):i\in \mathcal{R}\}$ is equal to the intersection of the images of $\gamma_{\max},\gamma_{\min}$. Indeed, from the fact that $\gamma_{\max},\gamma_{\min}$ are both geodesics started from $\varnothing$, $\gamma_{\max}(i)=\gamma_{\min}(j)$ automatically implies $i=j$. Recall that a (discrete) regenerative set is a random set of the form $\mathcal{G}=\{G_0+G_1+\ldots+G_n,n\geq 0\}$, where $G_1,G_2,\ldots$ are i.i.d.\ random variables with values in $\mathbb{N}$, independent of $G_0$. It is called aperiodic if the greatest common divisor of the support of the law of $G_1$ is $1$. If $\mathcal{G}$ is a regenerative set, then by the renewal theorem its asymptotic frequency exists and is given by
$$|\mathcal{G}|=\lim_{n\geq 1}\frac{\#\mathcal{G}\cap\{1,\ldots,n\}}{n}=\frac{1}{\mathbb{E}[G_1]}\, ,$$ and if $\mathcal{G}$ is aperiodic then $|\mathcal{G}|=\lim_{i\to\infty}\mathbb{P}(i\in \mathcal{G})$

\begin{lemma}\label{sec:proof-theorem-refcut} The set $\mathcal{R}$ is a (discrete) regenerative set, and $|\mathcal{R}|>0$ a.s.  \end{lemma}

\proof
We first note that the maximal geodesic only visits vertices in the trees $A_i,i\geq 0$, by (iii) in Lemma \ref{sec:maxim-minim-geod-1}.  More precisely, for every $i\geq 0$, $\gamma_{\max}(i)$ belongs to the tree $A_{j(i)}$ where $j(i)$ is the first index $j\leq i$ such that $\min\{\ell(v):v\in A_j\}=-i$.  Here we make a slight abuse of language by viewing $A_i$ as a (labeled) subgraph of $(T_\infty,\ell)$ rather than a tree in its own right.

In particular, we deduce that $\gamma_{\max}(i)=\gamma_{\min}(i)$ if and only if $j(i)=i$, that is $\inf\{\ell(v):v\in A_0\cup\ldots \cup A_{i-1}=-i+1$.
Otherwise said, letting $\Delta_i=\max\{-i-\ell(v):v\in A_i\}\geq 0$, then
$\max_{0\leq j\leq i-1}(j+\Delta_j)=i-1$, and this gives
$$ \mathcal{R}=\mathbb{Z}_+\setminus \bigcup_{j\geq 0}(j,j+\Delta_j]\, .$$
Yet otherwise said, $\mathcal{R}$ has same distribution as the set $\{G_0+G_1+\ldots+G_n:n\geq 0\}$ where $G_0=0$ and the variables $G_1,G_2,\ldots$ are i.i.d.\ with $G_1=\inf\{i>0:\max\{j+\Delta_j:0\leq j<i\}<i\}$.  Now, by Lemma \ref{indpt}, the random variables $(\Delta_i,i\geq 0)$ are independent and distributed as $\max \{-\ell(v):v\in \tau\}$ where $(\tau,\ell)$ has law $\rho_0$. By symmetry, this is the same as the law of $\max\{\ell(v):v\in \tau\}$, still under $\rho_0$. We will use the following lemma whose proof is postponed to Section \ref{scaling-limits}:
\begin{lemma}\label{sec:proof-theorem-refcut-1}
  It holds that, as $m\to \infty$,
 \begin{eqnarray*}\mathbb{P}(\Delta_0\geq m)&\sim& \frac{2}{m^2}\, , \end{eqnarray*}
the symbol $\sim$ meaning that the quotient of both sides tends to $1$. 
\end{lemma}

\noindent We note that \begin{eqnarray*}
  \mathbb{P}(i\in \mathcal{R})&=&\mathbb{P}(\max\{j+\Delta_j:0\leq j\leq i-1\}<i)\\
  &=&\prod_{j=0}^{i-1}\big(1-\mathbb{P}(\Delta_0\geq i-j)\big)\\
  &=&\prod_{j=1}^i\big(1-\mathbb{P}(\Delta_0\geq j)\big)\, , \end{eqnarray*} which as $i\to \infty$ decreases to a strictly positive limit by Lemma \ref{sec:proof-theorem-refcut-1} and the (obvious) fact that $\mathbb{P}(\Delta_0=0)>0$. This concludes the proof of Lemma \ref{sec:proof-theorem-refcut}.  
\endproof

Note that if $\mathcal{G}$ is a given infinite subset of $\mathbb{Z}_+$, then Lemma \ref{sec:proof-theorem-refcut} shows that $\mathcal{R}\cap \mathcal{G}$ is infinite almost surely, because 
$$ \mathbb{P}(\#\mathcal{R}\cap \mathcal{G}=\infty)\geq \lim_{i\to \infty,i\in \mathcal{G}}
\mathbb{P}(i\in \mathcal{R})>0\, ,$$ so the probability that $\#\mathcal{R}\cap \mathcal{G}=\infty$ has to be $1$ by the Hewitt-Savage $0$-$1$ law. 

\subsubsection{Leaving the spine} Our next step towards Theorem \ref{cut-points} is to prove that $\gamma_{\max}$ and $\gamma_{\min}$ eventually leave the spine for ever. We begin with the case of the maximal geodesic. 
Let 
 \begin{eqnarray*}\mathcal{R}'&=&\{i\geq 0:\gamma_{\max}(i)\in \{\mathrm{S}(j),j\geq 0\}\}  \end{eqnarray*}
be the set of times where $\gamma_{\max}$ hits the spine. 

\begin{lemma}
  \label{sec:proof-theorem-refcut-2}
  Almost surely, the set $\mathcal{R}'$ is finite. 
\end{lemma}

\proof
Let $L_0,L_1,L_2,\ldots$ be the subtrees to the left of the spine of $T_\infty$. Recall the notation $ \ell( \mathrm{S}(i)) = X_{i}$. Then note that the $i$-th vertex $\mathrm{S}(i)$ of the spine is on $\gamma_{\max}$ if and only if
$$\min\{\ell(v):v\in L_j\}> \ell(\mathrm{S}(i))\, ,\qquad \mbox{ for every }j\in \{0,1,\ldots,i-1\}\, .$$
Now,  for $i\geq 0$ let $\sigma_i=\inf\{n\geq 0:X_n=-i\}$, and let $L'_i$ be the forest made of the trees $L_{\sigma_i+j},0\leq j< \sigma_{i+1}-\sigma_i$. Then, if we let
$$\Delta'_i=\max_{v\in L'_i}-\ell(v)-i\, ,\qquad i\geq 0\, ,$$
then we obtain that
$$\mathcal{R}'=\mathbb{Z}_+\setminus
\bigcup_{j\geq 0}(j,j+\Delta'_j]\, .$$ Furthermore the $\Delta'_i$ are independent and identically  distributed. At this point the situation is similar to the setting of Lemma \ref{sec:proof-theorem-refcut} and we also rely on an external lemma whose proof is postponed to Section \ref{scaling-limits}:

\begin{lemma} \label{sec:proof-theorem-refcut-3}
As $m\to \infty$, it holds that 
 \begin{eqnarray*}\mathbb{P}(\Delta'_0\geq m)&\sim&\frac{2}{m}\, . \end{eqnarray*}
\end{lemma}

\noindent By the very same argument as for the computation of $\mathbb{P}(i\in \mathcal{R})$ we have using Lemma \ref{sec:proof-theorem-refcut-3}\begin{eqnarray*}
  \mathbb{P}(i\in \mathcal{R}')&=&\prod_{j=1}^i(1-\mathbb{P}(\Delta'_0\geq j))\\
  &=&\exp\Big(-\sum_{j=1}^i-\log(1-\mathbb{P}(\Delta_0'\geq j))\Big)\\
  &=&\exp\Big(-2\log m(1+o(1))\Big) \\&=&m^{-2+o(1)}\, .  \end{eqnarray*} In particular, $\mathbb{P}(i\in \mathcal{R}')$ is summable, so that, by the Borel-Cantelli Lemma, $i\notin\mathcal{R}'$ for every $i$ large enough, as wanted.  
\endproof

We deduce the following key property of the minimal geodesic. 

\begin{proposition}
   \label{sec:maxim-minim-geod}
Almost-surely, the minimal geodesic $\gamma_{\min}$ hits the spine a finite number of times: 
$$\#\big\{i\geq 0:\gamma_{\min}(i)\in \{\mathrm{S}(j),j\geq 0\}\big\}<\infty\, .$$
\end{proposition}

 \proof
   Assume by contradiction that $\gamma_{\min}$ hits the spine infinitely often with positive probability, and let $\mathcal{G}$ be the set of such intersections. Then $\mathcal{G}$ is measurable with respect to the subtree $\mathsf{C}$ obtained by chopping the trees $A_0,A_1,A_2,\ldots$ off $\theta$. By Lemma \ref{indpt}, we obtain that $\mathcal{G}$ is independent of $(A_i,i\geq 0)$, whence $\mathcal{G}$ is independent of the set $\mathcal{R}$ of intersections times of $\gamma_{\max}$ with $\gamma_{\min}$. Conditionally given $\mathcal{G}$, in restriction to the event that the latter is infinite, we then conclude from the  discussion after Lemma \ref{sec:proof-theorem-refcut} that $\mathcal{R}\cap \mathcal{G}$ is infinite almost-surely. But this is in contradiction with Lemma \ref{sec:proof-theorem-refcut-2}.
 \endproof

\subsubsection{Proof of Theorem \ref{cut-points}}

\proof[Proof of Theorem \ref{cut-points}] 
Let
 \begin{eqnarray*}i_*&=&\max\big\{j\geq 0:\gamma_{\min}(j)\in \{\mathrm{S}(0),\mathrm{S}(1),\ldots\}\big\}  \end{eqnarray*}
be the last time when $\gamma_{\min}$ hits the spine, 
and 
$$c_*=\min\{c:\mathcal{V}(c)=\gamma_{\min}(i_*)\}\, ,$$
the minimal corner incident to $\gamma_{\min}(i_*)$. Note that $c_*\leq c_0$, that is, $c_*$ lies on the right-hand side of $T_\infty$. Let $\gamma$ be a proper geodesic emanating from $\varnothing$. Recall from the proof of (iii) in Lemma \ref{sec:maxim-minim-geod-1} that $\gamma(i)$ belongs to the spine if and only if $\gamma(i)=\gamma_{\min}(i)$ and $\gamma_{\min}(i)$ belongs to the spine. Moreover,
after each visit to the spine, the corners used by $\gamma$ increase until the next visit to the spine. It follows that $c^{\gamma}(i)\geq c_*$ for every $i\geq 0$. In particular, $\gamma(i)$ is incident to the left-hand side of $T_\infty$ for every $i>i_0$, where $i_0=-\min\{\ell(c):c_*\leq c\leq c_0\}$ (the latter depends only on $(T_\infty,\ell)$ and not on the choice of $\gamma$). 

Now by (i) in Lemma \ref{sec:maxim-minim-geod-1}, we deduce that 
$c^{\gamma_{\max}}(i)\leq c^{\gamma}(i)\leq c^{\gamma_{\min}}(i)$ for every $i> i_0$. In particular, for every $i > i_{0}+1$ such that $i\in \mathcal{R}$, as defined around Lemma \ref{sec:proof-theorem-refcut}, we have $\gamma(i)=\gamma_{\max}(i)=\gamma_{\min}(i)$ by the discussion after the proof of Lemma \ref{sec:maxim-minim-geod-1}. Letting $p_1,p_2,\ldots$
be the ordered list of points of $\{\gamma_{\max}(i):i\in \mathcal{R},i>i_0\}$, which is infinite by Lemma \ref{sec:proof-theorem-refcut}, we thus see that every proper geodesic $\gamma$ has to visit the points $p_1,p_2,\ldots$
This concludes the proof of Theorem \ref{cut-points}. 
\endproof

\begin{figure}[htbp]
\begin{center} 
\includegraphics[width=16.5cm]{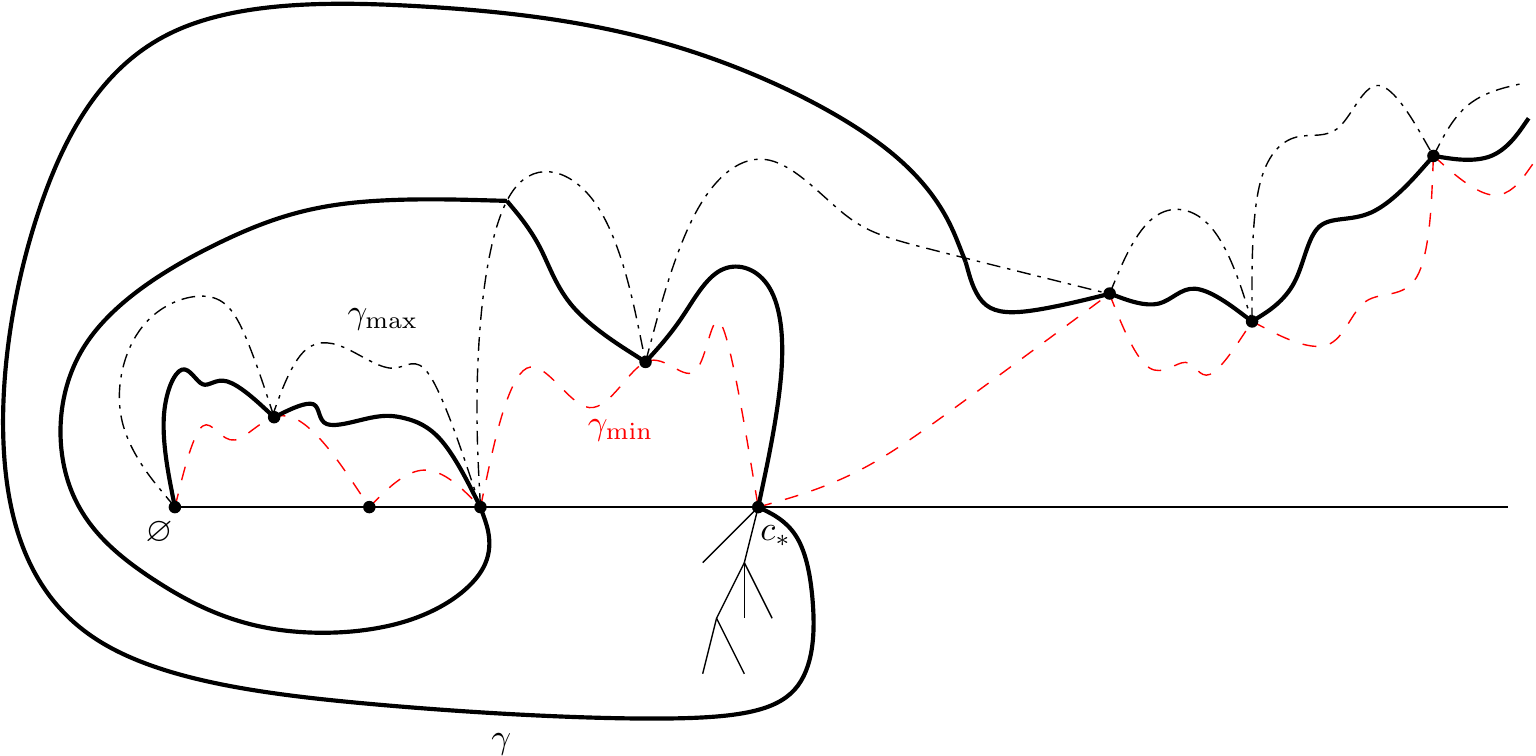}
 \end{center}
 \caption{Illustration of the proof of Theorem \ref{cut-points}: A proper geodesic can ``swirl'' around $T_\infty$ by traversing the spine, but is eventually trapped  between the minimal and maximal geodesics}
\label{escargot} 
 \end{figure}

\subsection{End of the proof of Theorem \ref{sec:unif-infin-label-1}}

\begin{lemma}\label{sec:proof-theor-refs}
Almost surely, the function $z\mapsto
\op{d}^{Q_\infty}_{\mathrm{gr}}(z,e^*_-)-\op{d}^{Q_\infty}_{\mathrm{gr}}(z,e^*_+)$
from $V(Q_\infty)$ to $\{-1,1\}$ is almost constant., $i.e.$ is constant except for finitely many $z \in Q_\infty$. 
\end{lemma}

\proof This statement is a property of the UIPQ, but for the purposes
of the proof, we will assume that $Q_\infty$ is constructed from a
tree $\theta=(T_{\infty}, \ell)$ with law $\mu$ and an independent parameter $\eta$
with $\mathcal{B}(1/2)$ distribution, by applying the Schaeffer
correspondence $\Phi$. This allows to specify the class of proper geodesic
rays among all geodesic rays emanating from $\varnothing$.

First, let us assume that $\varnothing=e^*_-$, meaning that
$\eta=0$. Let $\gamma_{\max}$ be the maximal geodesic so that
$\gamma_{\max}(0)=\varnothing=e^*_-$, and $\gamma_{\max}(1)=e^*_+$.  It is
also a proper geodesic ray, so that $\ell(\gamma_{\max}(i))=-i$ for
every $i \geq 0$.

Note that if $\gamma$ is a geodesic from $\varnothing$ to
$\gamma_{\max}(i)$ for some $i\geq 0$, then necessarily
$\ell(\gamma(j))=-j$ for every $j\in \{0,1,2,\ldots,i\}$, the reason
being that the labels of two neighboring vertices in $Q_{\infty}$ differ by at most
$1$.

Now let $\Gamma$ be the distinguished geodesic ray which starts from
$e^*_-=\varnothing$ constructed from $Q_\infty$ by first recovering the Chassaing-Durhuus tree $(\tau,\ell)= \overline{\Phi}{}^{-1}(Q_\infty)$ and then constructing $\Gamma$ as we did just before Lemma \ref{sec:confluent-geodesics-3}, and let $R\geq 0$. Applying Lemma 
\ref{sec:confluent-geodesics-3}, we obtain the existence of $R'\geq R$
such that the vertex $\gamma_{\max}(R'+1)$, which does not belong to
$B_{\carte,R'}(Q_\infty)$, can be linked to $\varnothing$ by a
geodesic $\gamma$ such that $\gamma(i)=\Gamma(i)$ for $i\in
\{0,1,\ldots,R\}$. Since
$\ell(\gamma(R'+1))=\ell(\gamma_{\max}(R'+1))=-(R'+1)$, we deduce from
the above discussion that $\ell(\gamma(i))=-i$ for every $i\in
\{0,1,\ldots,R'\}$, so in particular, $\ell(\Gamma(i))=-i$ for every
$i\in \{0,1,\ldots,R\}$. Since $R$ was arbitrary, we deduce that the
distinguished geodesic $\Gamma$ is proper.

By Theorem \ref{cut-points}, we get that $\Gamma$ and $\gamma_{\max}$
meet infinitely often. In particular, for every $\alpha\in (0,1)$, we
can find $R=R(\alpha)$ such that with probability at least $1-\alpha$,
there exists $I\in \{1,2,\ldots,R\}$ such that
$\Gamma(I)=\gamma_{\max}(I)$. From now on we argue on this event.
Applying Lemma \ref{sec:confluent-geodesics-3} again, we can
find $R'$ such that for every $z\in V(Q_\infty)\setminus
B_{\carte,R'}(Q_\infty)$, one can link $\varnothing$ to $z$ by a
geodesic $\gamma$ whose $R$ first steps coincide with those of
$\Gamma$. But since $\Gamma(I)=\gamma_{\max}(I)$, we can replace the
first $I$ steps of $\gamma$ by those of $\gamma_{\max}$, and obtain a
new geodesic from $\varnothing$ to $z$, whose first step goes from
$e^*_-$ to $e^*_+$. Since this holds for any $z$ at distance at least
$R'+1$ from $e^*_-$, we obtain that
$\op{d}^{Q_\infty}_{\mathrm{gr}}(z,e^*_-)-\op{d}^{Q_\infty}_{\mathrm{gr}}(z,e^*_+)=1$
for every $z$ at distance at least $R'+1$ from $e^*_-$. Since $\alpha$
was arbitrary, we obtain the desired result in the case $\eta=0$.

To treat the case $\eta=1$, we use the obvious fact that if
$\overleftarrow{Q}_\infty$ is the same quadrangulation as $Q_\infty$, but
where the root edge has the reverse orientation, then
$\overleftarrow{Q}_\infty$ has the same distribution as $Q_\infty$. Moreover,
$\overleftarrow{Q}_\infty=\Phi(\theta,1-\eta)$ so on the event
$\{\eta=1\}$ we are back to the situation $\eta=0$ by arguing on
$\overleftarrow{Q}_\infty$ instead of $Q_\infty$.  \endproof

From this, it is easy to prove (\ref{eq:4}), which will complete the
proof of Theorem \ref{sec:unif-infin-label-1}. Indeed, if $x$ and $y$ are neighboring vertices in $Q_{\infty}$ we can pick an edge $e$ such that $e_{-}=x$ and $e_{+}=y$. By Proposition  \ref{unimodular} below, the quadrangulation $Q_{\infty}^{(e)}$ re-rooted at $e$ has the same almost sure properties as $Q_{\infty}$. In particular, almost surely the function $z\mapsto
d(x,z)-d(y,z)$ is almost constant. But by reasoning on every step of a chain from
$x$ to $y$, the same holds for any $x,y \in Q_{\infty}$.  This constant has to be $\ell(x)-\ell(y)$. Indeed let us consider $\gamma_x$ and $\gamma_y$ two maximal geodesics emanating from a corner associated to $x$ resp.\,$y$. From properties of the Schaeffer construction, these two geodesics merge at some vertex $\gamma_x(i)= \gamma_y(i+ \ell(y)-\ell(x))$ for some $i \in\{0,1,2, \ldots\}$, and $\gamma_x(j)= \gamma_y(j+ \ell(y)-\ell(x))$ for every $j \geq i$. Hence 
$$ \lim_{z \to \infty} \left(\mathrm{d}_{ \mathrm{gr}}^{Q_\infty}(x,z)-\mathrm{d}_{ \mathrm{gr}}^{Q_\infty}(y,z)\right) =  \lim_{ {\begin{subarray}{c}   z\to \infty \\ z \in \gamma_x \cap \gamma_y \end{subarray}} } \left(\mathrm{d}_{ \mathrm{gr}}^{Q_\infty}(x,z)-\mathrm{d}_{ \mathrm{gr}}^{Q_\infty}(y,z) \right) =  \ell(x)- \ell(y).$$

\begin{corollary}\label{sec:coal-geod-rays}
Every geodesic ray emanating from $\varnothing$ is proper.
\end{corollary} 

\proof Let $\gamma$ be a geodesic ray and let $i_0 \geq 1$
fixed. Applying (\ref{eq:6}) we get 
\begin{eqnarray*}
\ell(\gamma(i_{0})) &=& \lim_{z\to
  \infty}\left(\op{d}_{\op{gr}}^{Q_{\infty}}(\gamma(i_{0}),z)-
\op{d}_{\op{gr}}^{Q_{\infty}}(\varnothing,z)\right)\\
&=&\lim_{i\to\infty}\left(\op{d}_{\op{gr}}^{Q_{\infty}}(\gamma(i_{0}),\gamma(i))-
\op{d}_{\op{gr}}^{Q_{\infty}}(\varnothing,\gamma(i))\right)\, .
\end{eqnarray*}
On the other hand, since $\gamma$ is a geodesic, for $i \geq i_{0}$ we
have $\op{d}_{\op{gr}}^{Q_{\infty}}(\gamma(i_{0}),\gamma(i)) =
i-i_{0},$ which implies that $\ell(\gamma(i_{0}))=-i_{0}.$ This allows
to conclude since $i_0$ was arbitrary.  \endproof

\section{Scaling limits for  $(T_{\infty}, \ell)$} \label{scaling-limits}
This section is devoted to the study of the scaling limits of the contour functions describing the tree $(T_\infty,\ell)$. It also contains the proof of  Lemmas \ref{sec:proof-theorem-refcut-1} and \ref{sec:proof-theorem-refcut-3} used during the proof of Theorem \ref{cut-points}.
\subsection{Contour functions}
\begin{figure}[!ht]
\begin{center}
\includegraphics[width=\textwidth]{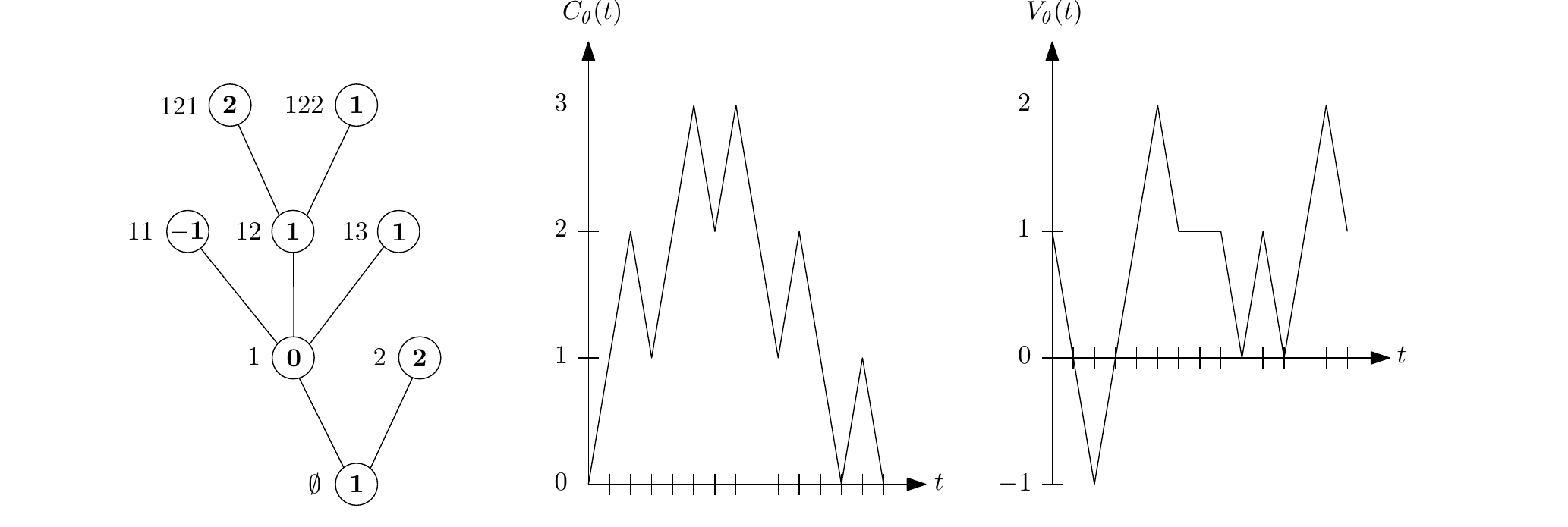}
\end{center}
\caption{A labeled tree $\theta$ and its pair of contour functions $\left( C_{\theta} ,V_{\theta} \right)$.}
\label{fig:contour}
\end{figure}

\paragraph{Coding of a tree.}
Let us recall a useful encoding of labeled trees by functions.  A finite labeled tree $\theta = (\tau, \ell)$ can be encoded by a pair $\left( C_{\theta},V_{\theta} \right)$, where $C_{\theta} = \left( C_{\theta} (t) \right)_{0 \leq t \leq 2|\theta|}$ is the contour function of $\tau$ and $V_{\theta} = \left( V_{\theta} (t) \right)_{0
  \leq t \leq 2|\theta|}$ is the labeling contour function of $\theta$. To define these contour functions, we let $(c_0,c_1\ldots,c_{2|\theta|-1})$ be the contour sequence of corners of $\tau$. Then $C_\theta(i)$ is the distance of $\mathcal{V}(c_i)$ to the root $\varnothing$ in $\tau$ for $0\leq i\leq 2|\theta|-1$, and we let $C_\theta(2|\theta|)=0$. Furthermore, we let $V_\theta(i)=\ell(c_i)$ for $0\leq i\leq 2|\theta|-1$, and then $V_\theta(2|\theta|)=\ell(c_0)$. Finally, we extend $C_\theta,V_\theta$ to continuous functions on $[0,2|\theta|]$ by linear interpolation between integer times (we will generally ignore these interpolations in what follows). 
See Figure \ref{fig:contour} for an example. A finite labeled tree is uniquely determined by its pair of contour functions. 

If  $\theta=(\tau,\ell)$ has law $\rho_l$ for some $l\in \mathbb{Z}$, then it is easy and well-known that $C_{\theta}$ is a simple random walk, stopped at time $\kappa_1-1$, where $\kappa_1$ the first hitting time of $-1$. For this reason, the process $C_\theta$ is sometimes extended to $[0,2|\theta|+1]$ by taking a final step of amplitude $-1$, this will help explain the construction of the process $C$ below. 
\paragraph{Coding of a forest.} If we now consider a sequence of trees $\theta_i=(\tau_i,\ell_i),i\geq 0$ respectively in $\ltrees{f}{l_i}$, 
 then we can concatenate the (extended) processes $C_{\theta_i}$ in a process
$$C(k)=C_{\theta_i}\Big(k-\sum_{j<i}(|2\theta_i|+1)\Big) - i \, ,\qquad \sum_{j<i}(2|\theta_j|+1)\leq k< \sum_{j\leq i}(2|\theta_j|+1)$$
We view $C$ as the contour function of an infinite forest $(\theta_0,\theta_1,\ldots)$, in which $C$ takes a $-1$ step at every newly visited tree. 
If we let $\underline{C}(i)=\inf_{0\leq j\leq i}C_j$, then the process $C_{\theta_i}$ can be recovered as
$$C_{\theta_i}(j)=i+C(\kappa_i+j) =C(\kappa_i+j)-\underline{C}(\kappa_i+j)\, ,\qquad 0\leq j< \kappa_{i+1}-\kappa_i=2|\theta_i|+1\, ,$$
where $\kappa_i=\inf\{n\geq 0:C(n)=-i\}$.  
We could also have chosen to take a $0$ or $+1$ steps at every newly visited tree: these contour functions are respectively given by $C- \underline{C}$ and $C-2\underline{C}$. We will see that it is easier to deal with $C$ but $C- 2 \underline{C}$ also plays a natural role, see the next paragraph.
\begin{figure}[!h]
 \begin{center}
 \includegraphics[width=14cm]{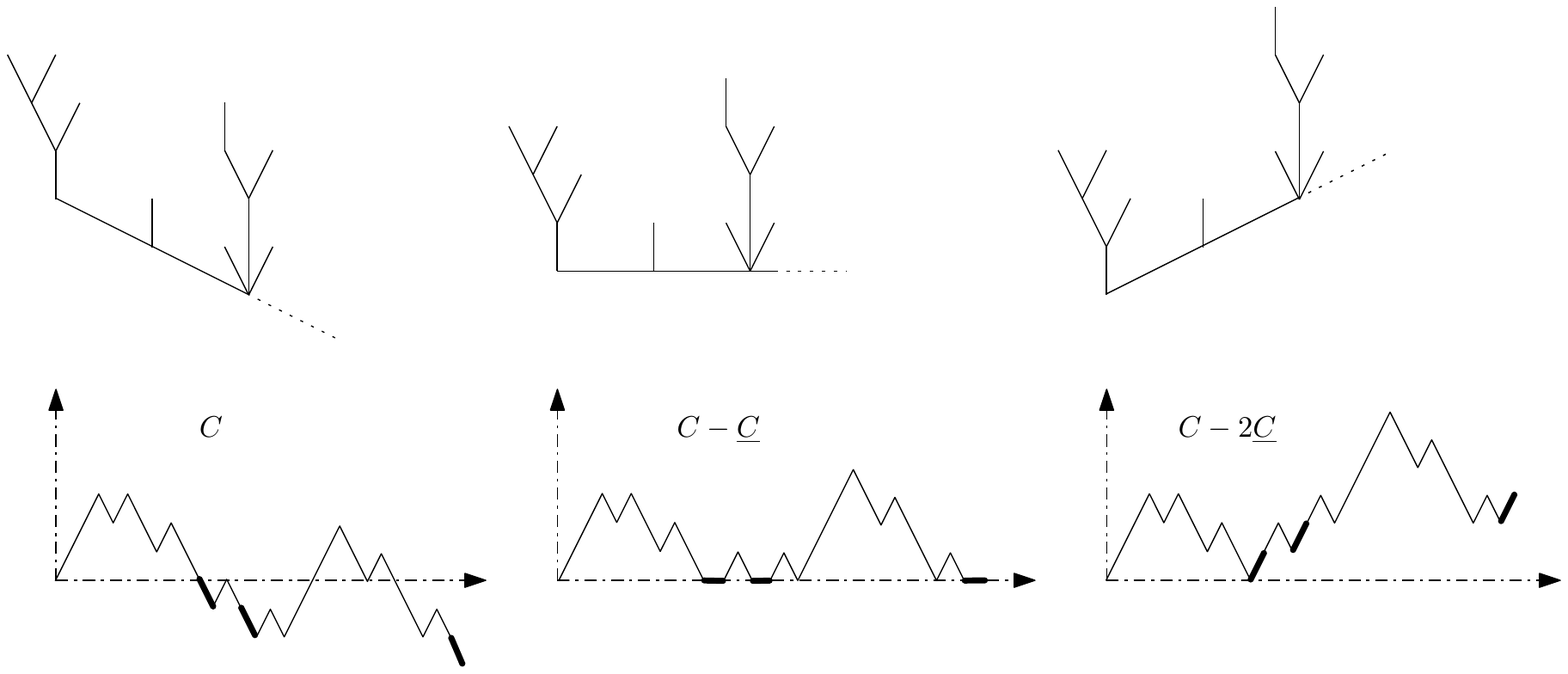}
 \caption{Different codings of a forest}
 \end{center}
 \end{figure}
 
As for the labeling function $V$, it is defined accordingly by concatenation of the processes $V_{\theta_i},i\geq 0$, so for $i\geq 0$, 
$$V(k+\kappa_i)=V_{\theta_i}(k)\, ,\qquad 0\leq k< \kappa_{i+1}-\kappa_i\,.$$

\subsection{Scaling limits of contour functions} \label{sec:scaling}

To describe the left part of the tree $T_{\infty}$ one thus would like to understand the contour functions associated to the sequence of trees $L_{0}, L_{1}, L_{2}, \ldots$. It turns out that it is easier to first deal with the contour and label processes of the labeled forest $L_0,L_1,\ldots$, in which we subtracted the label $X_i$ of the root of $L_i$ to the labels of all vertices of $L_i$. These sequences are denoted by $(C,V)$. Indeed, after doing this operation, the relabeled trees we obtain form an i.i.d.\ sequence with law $\rho_0$. In particular the process $C$ is a standard simple random walk with unit jumps.  The ``true'' label $\ell(\mathcal{V}(c_n))$ of the $n$-th vertex visited in the contour order of $T_\infty$ is then given by the formula \begin{equation}
  \label{eq:8}
  V'(n)=X_{-\underline{C}(n)}+V(n)\, ,
\end{equation}
because $-\underline{C}(n)$ is the index of the tree to which $\mathcal{V}(c_n)$ belongs. Obviously, the contour function $C$ is unchanged by this operation. Note that the labels on the spine $X$ is independent of $(C,V)$. We will introduce the scaling limits of $X$ and $(C,V)$ separatety. First, the Donsker invariance principle implies that 
  \begin{eqnarray}\left(\frac{1}{m}\sqrt{\frac{3}{2}}X_{m^2t}\right)_{t\geq 0}&\build\longrightarrow_{m\to\infty}^{(d)}& (\beta_t)_{t\geq 0}\, ,   \label{donskerfacile}\end{eqnarray}
where $\beta$ is a standard Brownian motion. On the other hand, let $(B_s,s\geq 0)$ be a standard Brownian motion, and $\underline{B}_s=\inf_{0\leq u\leq s}B_u$ be its infimum process. For $s \geq 0$, we let $\hat{B}_s=B_s-\underline{B}_s$, so that $\hat{B}=(\hat{B}_s,s\geq 0)$ is a standard reflected Brownian motion with local time at zero given by $(-\underline{B}_s,s\geq 0)$, by a famous theorem due to Lévy.  Now, conditionally given $\hat{B}$, let $(Z_s,s\geq 0)$
be a centered Gaussian process whose covariance function is given by
 \begin{eqnarray*} E[(Z_{s}-Z_{t})^2] &=& 
\hat{B}_s+\hat{B}_t-2\inf_{u\in [s\wedge t,s\vee t]}\hat{B}_u\, .
\end{eqnarray*}
The process $Z$ has a continuous modification \cite{LeG99} and this is the one we will deal with henceforth.  The pair $(\hat{B},Z)$ is called the Brownian snake (or sometimes the head of the Brownian snake) driven by the reflected Brownian motion $\hat{B}$. The reader can refer to the monograph \cite{LeG99} for a detailed account of the Brownian
snake. Then we have the joint convergence in distribution for the uniform norm over every compact interval:
 \begin{eqnarray}
\label{eq:7}
\left(\frac{1}{m^2}C(m^4t),\frac{1}{m^2}\underline{C}(m^4t),\frac{1}{m}\sqrt{\frac{3}{2}}V(m^4t)\right)_{t\geq 0}
&\build\longrightarrow_{n\to\infty}^{(d)}& (B_t,\underline{B}_t,Z_t)_{t\geq 0}\,,
 \end{eqnarray}
 this convergence holds jointly with \eqref{donskerfacile} and $(\beta_{t})$ is independent of the triplet $(B_{t},\underline{B}_{t},Z_{t})$. 
See for instance \cite[Theorem 3]{Mie08c} for a similar statement, from which the present one can be deduced easily. One difference is that \cite{Mie08c} deals with the so-called {\em height process} of the trees rather than the contour process, but the convergence of the latter is indeed a consequence of the convergence of height processes, as discussed in \cite[Section 2.4]{DLG02}. Convergences \eqref{donskerfacile} and \eqref{eq:7} entail that the contour functions $(C,V')$ of the forest $L_0,L_1, \cdots$ of the left part of $T_\infty$ admit the following scaling limit 
  \begin{eqnarray} \label{finalscaling} \left(\frac{1}{m^2}C(m^4t),\frac{1}{m}\sqrt{\frac{3}{2}}V'(m^4t)\right)_{t\geq 0}  & \xrightarrow[m\to\infty]{(d)} &(B_t,\beta_{-\underline{B}_t}+Z_t)_{t \geq 0}. \end{eqnarray}
To deal with the right part of $T_\infty$, we just remark that if $\tilde{C},\tilde{V}$ and $\tilde{V'}$ are defined from the forest $R_0, R_1, \cdots$ in the same way as $C,V$ and $V'$ were defined from $L_0, L_1, \ldots$ then we have the analogous of \eqref{eq:8} 
 \begin{eqnarray*} \tilde{V}'(n) &=&X_{-\underline{\tilde{C}}(n)}+\tilde{V}(n),  \end{eqnarray*} where $X$ is still the random walk of the labels on the spine. Remark that $(C,V)$ and  $(\tilde{C},\tilde{V})$ are independent and identically distributed and also independent of $X$. Hence the convergence \eqref{eq:7} also holds for $(\tilde{C},\tilde{V})$ namely jointly with \eqref{donskerfacile} and \eqref{eq:7} we have the convergence in distribution   \begin{eqnarray}
\label{eq:7bis}
\left(\frac{1}{m^2}\tilde{C}(m^4t),\frac{1}{m^2}\underline{\tilde{C}}(m^4t),\frac{1}{m}\sqrt{\frac{3}{2}}\tilde{V}(m^4t)\right)_{t\geq 0}
&\build\longrightarrow_{n\to\infty}^{(d)}& (\tilde{B}_t,\underline{\tilde{B}}_t,\tilde{Z}_t)_{t\geq 0}\,,
 \end{eqnarray}
 where $( \tilde{B},\underline{\tilde{B}}, \tilde{Z})$ is an independent copy of $(B, \underline{B},Z)$ also independent of $\beta$. Finally, the scaling limits of $\tilde{C}$ and $\tilde{V}'$ is given by an analogous formula as \eqref{finalscaling} after replacing $B,\underline{B}$ and $Z$ by their tilde-analogs. 
 
  \subsection{The continuous tree $ (\mathcal{T}_{\infty}, \mathcal{Z})$ }Let us give a slightly different point of view on these scaling limit results for the contour functions $(C,V)$ and $(\tilde{C},\tilde{V})$. The results of this paragraph will not be used in the sequel, we thus leave the details to the reader. 
  
  First of all, we remark that by a famous theorem of Pitman (see \cite[Chapter VI]{RY99}) the processes $B - 2 \underline{B}$ and $\tilde{B}-2\underline{\tilde{B}}$ are two independent three-dimensional Bessel processes.
These two processes thus give the scaling limits of the two contour functions of the left and right part of $T_{\infty}$ in which we make a $+1$ steps at the end of every visited tree. We now let 
  \begin{eqnarray*} \mathfrak{R}(t) &=& \left\{ \begin{array}{cl}
B_t-2\underline{B}_t
& \mbox{if } t \geq 0,\\
  \tilde{B}_{-t}-2\underline{\tilde{B}}_{-t} & \mbox{if } t \leq 0,
  \end{array} \right.
  \end{eqnarray*}
and for every $s,t \in \mathbb{R}$, we define
   \begin{eqnarray*} \overline{st} &=& \left\{ \begin{array}{ll}
   [s\wedge t, s \vee t] \hspace{2cm} & \mbox{if } st\geq 0,\\
   (-\infty, s \wedge t] \cup [s \vee t, \infty) & \mbox{if }st  <0. \end{array}  \right.\end{eqnarray*}
 Finally we define a pseudo-distance on $ \mathbb{R}$ by letting  \begin{eqnarray*}
   \mathrm{d}_{ \mathfrak{R}}(s,t) &=& \mathfrak{R}_{s}+ \mathfrak{R}_{t}-2 \inf_{ r\in \overline{st}} \mathfrak{R}_{r}.  \end{eqnarray*} The quotient space $ \mathcal{T}_{\infty} = \mathbb{R} / (\mathrm{d}_{ \mathfrak{R}}=0)$ equipped with the quotient distance $ \mathrm{d}_{ \mathfrak{R}}$ is an infinite tree that is the scaling limit of the tree $T_{\infty}$, that is $\lambda \cdot T_{\infty} \to \mathcal{T}_{\infty}$ in distribution for the Gromov-Hausdorff distance \cite{BBI01} as $\lambda \to 0$. 
   
   Furthermore, conditionally on $ \mathfrak{R}$, we consider a real-valued centered Gaussian process $ \mathcal{Z}_{t}$ indexed by $ \mathbb{R}$ whose covariance function is prescribed by  \begin{eqnarray*} E[( \mathcal{Z}_{s}- \mathcal{Z}_{t})^2] &=& \mathrm{d}_{ \mathfrak{R}}(s,t).  \end{eqnarray*}
   Similarly to the case of $Z$, the process $ \mathcal{Z}$ has a continuous modification that we consider from now on. In words, the process $ \mathcal{Z}$ can be interpreted as the Brownian motion indexed by the infinite tree $ \mathcal{T}_{\infty}$. We claim that conditionally on $C$ and on $ \tilde{C}$ we have the equality in distribution 
    \begin{eqnarray*} \Big( \mathbf{1}_{t \geq 0} \big(Z_{t} + \beta_{-\underline{B}_{t}}\big) + \mathbf{1}_{t \leq 0} \big(\tilde{Z}_{-t} + \beta_{-\underline{ \tilde{B}}_{-t}}\big) \Big)& \overset{(d)}{=} & (\mathcal{Z}_{t})_{t \in \mathbb{R}},  \end{eqnarray*}
    which can be check by looking at the covariance functions of these Gaussian processes. To conclude, this interpretation enables us to fully understand that the labeled tree $T_{\infty}$ converges, in the scaling limit, towards a non-compact random tree $ \mathcal{T}_{\infty}$ encoded by a pair of independent Bessel processes with the Brownian labeling $ \mathcal{Z}$. This object should play a crucial role to describe the scaling limit of the UIPQ in the Gromov-Hausdorff sense. We plan to study this in future work.

\subsection{Proofs of Lemmas \ref{sec:proof-theorem-refcut-1} and \ref{sec:proof-theorem-refcut-3}}

In this section we proceed to the proof of the lemmas used during the proof of Theorem \ref{cut-points}. Although these lemmas seem combinatorial in nature, the  constant appearing in the equivalents are obtained by using the continuous scaling limits of the contour processes of the tree $T_\infty$. \medskip

Recall that $\Delta_0$ is distributed as $\max\{\ell(v) : v \in \tau\}$ under $\rho_0$. \medskip 

\noindent{\bf Lemma \ref{sec:proof-theorem-refcut-1}.}  \emph{
  It holds that, as $m\to \infty$,
 \begin{eqnarray*}\mathbb{P}(\Delta_0\geq m)&\sim& \frac{2}{m^2}. \end{eqnarray*}}

\proof Consider an infinite sequence $(\theta_i,\ell_i)_{i\geq 0}$ of independent labeled trees with law $\rho_0$. 
We let $(C,V)$ be the contour and label sequences of this forest, as defined in the last section. Recall the notation $\kappa_i=\inf\{n\geq 0:C(n)=-i\}$. The convergence \eqref{eq:7}, together with standard arguments relying on the fact that for a given level $-x$, the first hitting time $T_x=\inf\{t\geq 0:B_t<-x\}$ (the definition is chosen so that $T$ is the right-continuous inverse of the function $-\underline{B}$) is almost surely not a local minimum for the function $B$, entails that 
 \begin{eqnarray}
  \label{eq:12}
\frac{1}{m}\sqrt{\frac{3}{2}}\max\{V(i):0\leq i\leq \kappa_{m^2}\}&\build\longrightarrow_{m\to \infty}^{(d)} & \sup\{Z_t:0\leq t\leq T_1\}\, .
 \end{eqnarray}
By excursion theory for the Brownian snake \cite{LeG99}, it holds that
$$Z^{(a)}=(Z_{T_{a-}+s},0\leq s\leq T_{a}-T_{a-})\, ,\qquad a\geq 0$$
is a Poisson point process with intensity $2\mathbb{N}_0(\d Z)$, where $\mathbb{N}_0$ is the so-called Itô excursion measure of the Brownian snake, and with standard Brownian spatial displacements. In particular, it is known, by \cite[Lemma 2.1]{LGW06} and invariance in distribution of the process $Z$ under reflection, that
$$2\mathbb{N}_0(\sup Z>x)=\frac{3}{x^2}\, ,\qquad x>0\, .$$
Therefore, by standard properties of Poisson random measures, 
$$\mathbb{P}(\sup \{Z_t:0\leq t\leq T_1\}\leq x)=\mathbb{P}(\sup \{\sup Z^{(a)}:0\leq a\leq 1\}\leq x)=\exp(-3/x^2)\, .$$ 
Moreover, since the paths $(V(\kappa_i+j),0\leq j\leq \kappa_{i+1}-\kappa_i),i\geq 0$ are i.i.d.\ and $\max \{V(j):0\leq j\leq \kappa_1\}$ has same law as $\Delta_0$, we obtain that for every $x>0$, 
\begin{eqnarray*}
  \left(1-\mathbb{P}\Big(\Delta_0>mx\sqrt{\frac{2}{3}}\Big)\right)^{m^2+1}&=&
  \mathbb{P}\left(\frac{1}{m}\sqrt{\frac{3}{2}}\max \{V(i):0\leq i\leq \kappa_{m^2}\}\leq x\right)\\
  &\build\longrightarrow_{m\to\infty}^{}&\exp(-3/x^2)\, ,
\end{eqnarray*}
where in the last step we used \eqref{eq:12} and the fact that the random variable $\sup \{Z_t:0\leq t\leq T_1\}$ has a diffuse law.  This entails
$$m^2\mathbb{P}(\Delta_0>m)\to  2\, ,$$
concluding the proof of Lemma \ref{sec:proof-theorem-refcut-1}.  \endproof 

Next, using the notation of Section \ref{sec:scaling}, we let
$$\kappa'_i=\inf\{j\geq 0:X_{-\underline{C}(j)}=-i\}\, .$$ We recall the definition of the sequence $(\Delta'_i,i\geq 0)$ introduced in Lemma \ref{sec:proof-theorem-refcut-2} 
\begin{eqnarray*}
  \Delta'_i&=&\max\{-V'(\kappa_{i}'+j),0\leq j\leq \kappa'_{i+1}-\kappa'_i\}-i.\nonumber
   \end{eqnarray*}

\noindent {\bf Lemma  \ref{sec:proof-theorem-refcut-3}.} \emph{
As $m\to \infty$, it holds that 
 \begin{eqnarray*}\mathbb{P}(\Delta'_0\geq m)&\sim&\frac{2}{m}\, . \end{eqnarray*}}

\proof  
  The sequence $\Delta'_i$ admits the alternative representation
  \begin{eqnarray}
\Delta_i'&=& \max\{-V(\kappa_{i}'+j)-X_{-\underline{C}(\kappa'_i+j)}+\underline{X}_{-\underline{C}(\kappa'_i+j)},0\leq j\leq \kappa'_{i+1}-\kappa'_i\}\, ,\qquad i\geq 0\, ,\label{eq:11}
\end{eqnarray}
where $\underline{X}_i=\min\{X_j,0\leq j\leq i\}$, so that
$$\max_{0\leq i\leq m}\Delta'_i=\max\{-V(j)-X_{-\underline{C}(j)}+\underline{X}_{-\underline{C}(j)}:0\leq j\leq \kappa'_m\}\, .$$
Note also that $(\Delta'_i,i\geq 0)$ is an i.i.d.\ sequence. This can be proved by exploring the sequene $L_0, L_1,\ldots$ in a Markovian way, in the spirit of the proof of Lemma \ref{indpt}, and we leave this fact to the reader. 
Convergences \eqref{donskerfacile} and \eqref{eq:7} entail that 
$$\Big(\frac{1}{m}\sqrt{\frac{3}{2}}(V(m^4t)+X_{-\underline{C}(m^4t)}-\underline{X}_{-\underline{C}(m^4t)})\Big)_{t\geq 0}\build\longrightarrow_{m\to\infty}^{} (Z_t+\hat{\beta}_{-\underline{B}_t}\,t\geq 0)\, ,$$
where $\hat{\beta}_x=\beta_x-\underline{\beta}_x$.
Let $T'_x=\inf\{t\geq 0:\beta_{-\underline{B}_t}<-x\}=T_{\tau_x}$, where $\tau_x=\inf\{u\geq 0:\beta_u<-x\}$ is the right-continuous inverse of $-\underline{\beta}$. For a given $x\geq 0$, it is easy to check that the time $T'_x$ is not a local minimum for the process $(\beta_{-\underline{B}_t},t\geq 0)$.  Therefore, the previous convergence entails
\begin{equation}
  \label{eq:10}
  \Big(\frac{1}{m}\sqrt{\frac{3}{2}}(V'(m^4t) -\underline{X}_{-\underline{C}(m^4t)})\Big)_{0\leq t\leq \kappa'_m/m^4}\build\longrightarrow_{m\to\infty}^{} (Z_t+\hat{\beta}_{-\underline{B}_t})_{0\leq t\leq T'_{\sqrt{3/2}}}\,
\end{equation}
Let us fix $x\geq 0$. We claim that
\begin{equation}
  \label{eq:9}
  \mathbb{P}\Big(\sup\big\{Z_t-\hat{\beta}_{-\underline{B}_t}:0\leq t\leq T'_{\sqrt{3/2}}\big\}\leq x\Big)=\exp(-\sqrt{6}/x)\, .
\end{equation}
 To show this, recall the notation $Z^{(a)}$ used in the proof of Lemma \ref{sec:proof-theorem-refcut-1}. Then, \begin{eqnarray*}
  \lefteqn{  \sup\big\{Z_t-\hat{\beta}_{-\underline{B}_t}:0\leq t\leq T'_{\sqrt{3/2}}\big\}
  }\\
  &=&\sup\{\sup\{Z_{T_{a-}+s}:0\leq s\leq T_{a}-T_{a-}\}-\hat{\beta}_a:0\leq a\leq \tau_{\sqrt{3/2}}\}\\
  &=&\sup\{\sup Z^{(a)}-\hat{\beta}_a:0\leq a\leq \tau_{\sqrt{3/2}}\}\, .
\end{eqnarray*}
Since $(Z^{(a)},a\geq 0)$ is a homogeneous Poisson process, standard properties of Poisson measures entail that conditionally given $\beta$, 
\begin{eqnarray*}
 \mathbb{P}\Big(\sup\big\{Z_t-\hat{\beta}_{-\underline{B}_t}:0\leq t\leq T'_{\sqrt{3/2}}\big\}\leq x\, \Big|\, \beta\Big)&=&  \exp\left(-\int_0^{\tau_{\sqrt{3/2}}}\d a\, 2\mathbb{N}_{0}(\sup Z-\hat{\beta}_a\geq x)\right)\\
&=&\exp\left(-\int_0^{\tau_{\sqrt{3/2}}}\frac{\d a}{(x+\hat{\beta}_a)^2}\right)\\
&=&\exp\left(-\sum_{0<y\leq \sqrt{3/2}}\int_0^{\tau_y-\tau_{y-}}\frac{\d a}{(x+\beta^{(y)}_a)^2}\right)
\end{eqnarray*}
where
$$\beta^{(y)}_a=\hat{\beta}_{\tau_y+a}=y+\beta_{\tau_y+a}\, ,\qquad 0\leq a\leq \tau_y-\tau_{y-}$$
is the excursion of $\beta$ above its minimum at level $-y$. Taking expectations and using Itô's excursion theory, we obtain
$$ \mathbb{P}\Big(\sup\big\{Z_t-\hat{\beta}_{-\underline{B}_t}:0\leq t\leq T'_{\sqrt{3/2}}\big\}\leq x\Big)=\exp\Big(-\sqrt{\frac{3}{2}}\int 2n(\d e)\int_0^{\sigma(e)}\frac{\d a}{(x+e_a)^2}\Big)\, ,$$
where $n(\d e)$ is the Itô measure of positive excursions of Brownian motion (so that $2n$ is the excursion measure of the reflected Brownian motion $\hat{\beta}$), and $\sigma(e)$ is the lifetime of the generic excursion $e$. The Bismut decomposition of the measure $n$ finally gives
$$\int n(\d e)\int_0^{\sigma(e)}\frac{\d a}{(x+e_a)^2}=\int_0^\infty\frac{\d r}{(x+r)^2}=\frac{1}{x}\, ,$$
hence \eqref{eq:9}. Putting \eqref{eq:11}, \eqref{eq:10} and \eqref{eq:9} together, we deduce that 
\begin{eqnarray*}
  \Big(1-\mathbb{P}(\Delta'_0>xm\sqrt{\frac{2}{3}})\Big)^{m+1}&=&
  \mathbb{P}\Big(\frac{1}{m}\sqrt{\frac{3}{2}}\max\{\Delta'_i:0\leq i\leq m\}
  \leq x\Big)\\  &=&\mathbb{P}\Big(\frac{1}{m}\sqrt{\frac{3}{2}}\max\{-V'(j)+\underline{X}_{-\underline{C}(j)}:0\leq i\leq \kappa'(m)\}\leq x\Big)\\
  &\build\longrightarrow_{m\to\infty}^{}& \mathbb{P}\Big(\sup\big\{Z_t-\hat{\beta}_{-\underline{B}_t}:0\leq t\leq T'_{\sqrt{3/2}}\big\}\leq x\Big)=\exp(-\sqrt{6}/x)\, ,
\end{eqnarray*}
where at the penultimate step we used the fact that the law of $\sup\{Z_t-\hat{\beta}_{-\underline{B}_t}:0\leq t\leq T'_{\sqrt{3/2}}\}$ is diffuse, and the fact that $(B,-Z,-\beta)$ and $(B,Z,\beta)$ have the same distribution.
Therefore, taking $x=\sqrt{3/2}$, we get
$$\mathbb{P}(\Delta_0'>m)\sim \frac{2}{m}\, ,$$
concluding the proof of Lemma \ref{sec:proof-theorem-refcut-3}. \endproof

\section{Horoballs and points of escape to infinity}\label{sec:horoballs}

In the two remaining sections, we use the representation of the UIPQ given by Theorem \ref{sec:unif-infin-label-1} in order to deduce new results on this object.  First, inspired by the work of Krikun \cite{Kri05}, we study the length of the separating cycle around the origin of the UIPQ at ``height'' $-r$ seen from infinity.
In the following, we let $Q_\infty$ be a random variable with the law of the UIPQ, and we let $((T_\infty,\ell),\eta)$ be the labeled tree associated with $Q_\infty$ by the Schaeffer correspondence.  \bigskip

For every integer $l \in \mathbb{Z}$, we denote by $H_l=H_l(Q_\infty)=\left\{ v \in V(Q_\infty) : \ell(v) \leq  l\right\}$. In view of Theorem \ref{sec:unif-infin-label-1}, we can interpret $H_l$ as the {\em ball centered at infinity}, or {\em horoball}, with ``radius $\op{d}^{Q_\infty}_{\mathrm{gr}}(\varnothing,\infty)+l$'', where $\infty$ is a point at infinity in $Q_\infty$. Intuitively, the boundary of this set (a ``horosphere'') is made of several disjoint cycles, one of which separates $\varnothing$ from $\infty$. We are going to give asymptotic properties for the length of this cycle, the set of ``points to escape to $\infty$ at level $l$''. 

To this purpose, it is easier to work with a slightly modified graph $\hat{Q}_\infty$, which is obtained by adding to $Q_\infty$ the edges of $T_\infty$ given by the inverse Schaeffer construction of Section \ref{quadtotrees}. Therefore, in faces $f$ of $Q_\infty$, around which the four vertices in clockwise order have labels $l,l+1,l,l+1$, we add the diagonal between the vertices with label $l+1$ (such faces are called confluent faces in \cite{CS04}).  In faces with labels $l,l+1,l+2,l+1$, we just double the edge between the last two vertices. The map $\hat{Q}_\infty$ is is no longer a quadrangulation: Some of its faces are still squares, but others are triangles, and some others have degree $2$. Nonetheless, $Q_\infty$ and $\hat{Q}_\infty$ are very similar: They have the same vertex set, and it is easy to check that formula \eqref{eq:6} remains true in this context, i.e. 
\[ \ell(v) = \lim_{z \to \infty} \left({\op d}_{\mathrm{gr}}^{Q_\infty}(z, \varnothing) - \op{d}_{\mathrm{gr}}^{Q_\infty} (z,v) \right)= \lim_{z \to \infty} \left(\op{d}_{\mathrm{gr}}^{\hat{Q}_\infty} (z, \varnothing) - \op{d}_{\mathrm{gr}}^{\hat{Q}_\infty} (z,v)\right) \] for every vertex $v$, because we only added edges between vertices with the same label, and geodesic rays in $\hat{Q}_\infty$ never use the new edges. We are not going to need this fact in the sequel, so details are left to the reader.

Now the complement $\{v\in V(Q_\infty):\ell(v)> l \}$ of $H_l$ induces a subgraph of $\hat{Q}_\infty$, and if $l<0$ we let $F_l=F_l(Q_\infty)$ be the set of vertices in the connected component of this subgraph that contains $\varnothing$.  We let $\partial F_l$ be the set of vertices of $H_l$ that are connected to $F_l$ by an edge.
Yet otherwise said, $F_l$ is the set of vertices $v$ with $\ell(v)> l$, and which can be joined to $\varnothing$ by a path of $\hat{Q}_\infty$ along which labels are all strictly greater than $ l$, and $\partial F_l$ is the set of vertices with label $l$ that can be joined to $\varnothing$ by a path in $Q_\infty$ along which all labels are all strictly greater than $l$ except at the initial point.

Although it is not obvious at first sight, the set $F_l$ is almost-surely finite: The next statement implies that it is contained in the set of vertices of the subtrees $L_0,\ldots,L_{\sigma_{l}},R_0,\ldots,R_{\sigma_{l}}$ to the left and to the right of $T_\infty$, where $\sigma_{l}=\inf\{n\geq 0:\ell(\mathrm{S}(n))=-l\}$ is a.s.\ finite.
Recall that $\llbracket a,b\rrbracket$ is the path from $a$ to $b$ in $T_\infty$.

\begin{proposition}
\label{lem:exit}
Let $r>0$ and $v\in V(Q_\infty)$. Then 
\begin{enumerate}
\item 
$v$  belongs to $F_{-r}(Q_\infty)$ if and only if $\ell(v')>-r$ for every $v'\in \llbracket\varnothing, v\rrbracket$, and  
\item
$v$  belongs to $\partial F_{-r}(Q_\infty)$ if and only if $\ell(v)=-r$ and $\ell(v')>-r$ for every $v'\in \llbracket\varnothing, v\rrbracket\setminus \{v\}$. 
\end{enumerate}
\end{proposition}

\proof
  Let $v$ be such that $\ell(v')>-r$ for every $v'\in \llbracket \varnothing,v\rrbracket\setminus\{v\}$. The path $\llbracket \varnothing,v\rrbracket$ in $T_\infty$ is also a path in the augmented graph $\hat{Q}_\infty$, and it goes only through vertices outside $H_{-r}$, except maybe at $v$. So if furthermore $\ell(v)=-r$ we obtain that $v\in \partial F_{-r}$, otherwise $v\in F_{-r}$. 

Conversely, suppose that  there is a vertex
$v'$ on the path $[[\varnothing, v]]\setminus \{v\}$ with $\ell(v') \leq -r$. 
 We consider the two corners $c_{(1)}$ and $c_{(2)}$ that are respectively the smallest and largest corner incident to the vertex $v'$. We then construct the geodesics $\gamma_1$ and
$\gamma_2$ emanating from the corners $c_{(1)}$ and $c_{(2)}$ by taking consecutive successors $\mathcal{S}^{(i)}(c_{(1)}),i\geq 0$ and $\mathcal{S}^{(i)}(c_{(2)}),i\geq 0$ respectively. 

\begin{figure}[!h]
\begin{center}
\includegraphics[width=0.80\textwidth]{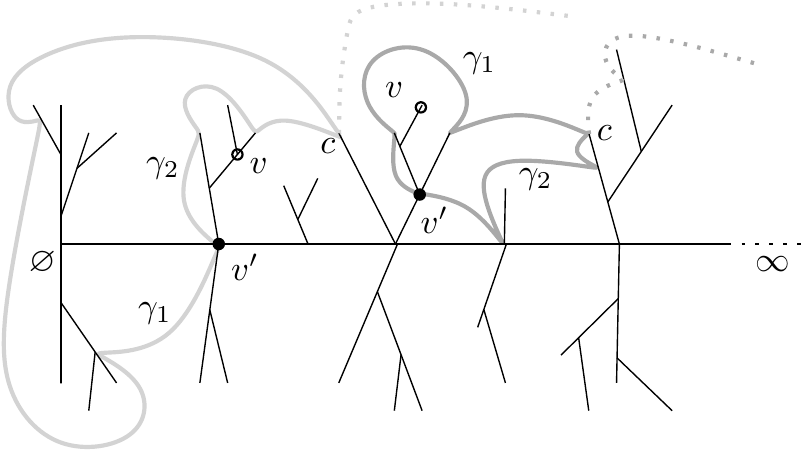}
\caption{Cycles separating vertices from $\varnothing$
composed by two geodesics.}
\label{cycleexit}
\end{center}
\end{figure}

 The geodesics $\gamma_1$ and $\gamma_2$ coalesce at the first corner $c$ following $c_{(2)}$ in the contour with label $\ell_{\min} - 1$, where $\ell_{\min}$ is the smallest label of the corners between $c_{(1)}$ and $c_{(2)}$.  The concatenation of the parts of $\gamma_1$ and $\gamma_2$ between $v'$ and $\mathcal{V}(c)$ induces a cycle ${\mathcal C}$ of the map $Q_\infty$, such that $v \notin \mathcal{C}$ and separates $v$ from $\varnothing$ in $Q_\infty$, see Figure \ref{cycleexit} for an example. Note that all the vertices of $ \mathcal{C}$ have labels less than or equal to $-r$. By the Jordan Theorem, any path in $\hat{Q}_\infty$ joining $\varnothing$ to $v$ crosses $\mathcal C$, and thus has a vertex, other than $v$, with label
 less than or equal to $-r$. It follows that $v$ does not belong to $\partial F_{-r}\cup F_{-r}$.
\endproof

We can use Proposition \ref{lem:exit} to derive asymptotics for the number of
vertices in $\partial F_{-r} (q)$:

\begin{theorem} \label{separatingcycles}
The sequence $(2|\partial F_{-r}(Q_\infty)|/r^2,r\geq 1)$ converges in distribution to a random variable $ W$ with Laplace transform
 \begin{eqnarray*}
E[e^{-\lambda W}]&=&
\left( \frac{1}{1 + \sqrt{\lambda}} \right) ^3\, .
 \end{eqnarray*}
\end{theorem}

\begin{remark} 
A similar problem has already been studied by Krikun. In  \cite{Kri05}, he considered the component $\partial\tilde{F}_r(Q_\infty)$ of the boundary of the ball $B_{\carte,r}(Q_\infty)$ that separates the root vertex from infinity. He showed that  $2|\partial\tilde{F}_r(Q_\infty)|/r^2$ converges in distribution as $r\to\infty$ to a standard Gamma random variable with parameter $3/2$.  In our setup, $|\partial F_{-r}(Q_\infty)|$ is, roughly speaking, the length of the boundary of the horoball that separates the root from infinity at level $-r$. It is not surprising that the typical size of this component should be of order $r^2$ as well, and obey a similar limiting result. However, we see that the limiting distribution differs from the Gamma law. In this respect, it is interesting to compare the tail distributions of these variables. Let $W,\tilde{W}$ have respective Laplace transforms
$$\phi(\lambda)=E[e^{-\lambda W}]=\Big(\frac{1}{1+\sqrt{\lambda}}\Big)^3\, ,\qquad 
\tilde{\phi}(\lambda)=E[e^{-\lambda \tilde{W}}]=\left(\frac{1}{1+\lambda}\right)^{3/2}\, ,$$ so that $\tilde{W}$ follows the Gamma distribution with parameter $3/2$.  As $\lambda\to \infty$, we have $\phi(\lambda)\sim \tilde{\phi}(\lambda)\sim \lambda^{-3/2}$, so that a Tauberian Theorem \cite[Theorem 1.7.1']{BGT89} entails that as $y\to 0^+$
$$\mathbb{P}(W\leq y)\sim \mathbb{P}(\tilde{W}\leq y)\sim \frac{4}{3\sqrt{\pi}}y^{3/2}\, ,$$
so these two distributions have a similar behavior at $0$. By contrast, $\tilde{W}$ has exponential tails, while we have 
$$\phi(\lambda)=1-3\sqrt{\lambda}+o(\sqrt{\lambda})\, , \quad \mbox{ when }\lambda \to 0 $$
which by applying Corollary 8.1.7 in \cite{BGT89} shows that
$$\mathbb{P}(W\geq y)\sim \frac{3}{\sqrt{\pi}}y^{-1/2}\, ,$$
as $y\to\infty$. Therefore, $W$ has a heavy tail. 
\end{remark}

\proof[Proof of Theorem \ref{separatingcycles}] Assume that $Q_{\infty}$ is obtained by the extended Schaeffer bijection, that is $Q_\infty=\Phi((T_{\infty},\ell),\eta)$ where $((T_{\infty},\ell),\eta)$ has law $\mu\otimes \mathcal{B}(1/2)$. The only vertex of $\partial F_{-r}(Q_{\infty})$ that belongs to the spine of $T_{\infty}$ is $ \mathrm{S}( \sigma_{r})$ where we recall that $\sigma_{r} = \inf\{ i \geq 0: \ell( \mathrm{S}(i))=-r\}$. Then, Proposition \ref{lem:exit} implies that \begin{equation} \label{eq:decompexit} |\partial F_{-r}(Q_\infty)| = 1 + \sum_{i < \sigma_{r}} \left| Y_{L_i}(-r) \right| + \sum_{i < \sigma_{r}} \left| Y_{R_i}(-r) \right| \end{equation} where, if $\theta$ is a labeled tree whose root label is strictly larger than $-r$, $Y_{\theta}(-r)$ is the set of all vertices of $\theta$ with label $-r$ and such that all their ancestors have a label strictly larger than $-r$.  Recalling that the trees $(L_i)$ and $(R_i)$ are independent conditionally given $(X_i)$, we get from \eqref{eq:decompexit} that \begin{align} \label{laplaceexit}
      E  &\left[ \exp - \lambda \frac{2|\partial F_{-r}(Q\infty)|}{r^2}\right] \notag\\
      &= \exp \left( - \frac{2\lambda}{r^2} \right) E \left[ \prod_{i < \sigma_{r}} E \left[\exp \left( - \frac{2\lambda}{r^2} \left| Y_{L_i}(-r) \right| \right) \middle| X_i \right] E \left[\exp \left( - \frac{2\lambda}{r^2} \left| Y_{R_i}(-r) \right| \right) \middle| X_i \right]
      \right] \notag\\
      & = \exp \left( - \frac{2\lambda}{r^2} \right) E \left[ \prod_{i < \sigma_{r}} \rho_{X_i} \left( \exp \left( - \frac{2\lambda}{r^2} \left| Y_{\theta}(-r) \right| \right) \right)^2
      \right] \notag\\
      & = \exp \left( - \frac{2\lambda}{r^2} \right) E \left[ \exp \left( 2 \sum_{i < \sigma_{r}} \log \rho_{X_i} \left( \exp \left( - \frac{2\lambda}{r^2} \left| Y_{\theta}(-r) \right| \right) \right) \right) \right].  \end{align}

\bigskip

To compute the right hand side of \eqref{laplaceexit}, we need to evaluate the
generating functions:
\[
f_{l,-r}(x) = \rho_{l} \left( x^{\left| Y_{\theta}(-r) \right|} \right), \,
0 \leq x \leq 1
\]
for $l \geq -r$, with the boundary condition $f_{-r,-r}(x)=x$. 
The measures $\rho_l$ being the laws of
Galton Watson trees with geometric$(1/2)$ offspring distribution and uniform labels, it is easy to derive the following
recursive relation for $l> -r$:
\begin{eqnarray*}
f_{l,-r}(x) &=& \frac{1}{2}\sum_{k = 0}^{\infty} \, \sum_{i_1+i_2+i_3 = k}
\binom{k}{i_1,i_2,i_3}
\left( \frac{1}{6} f_{l-1,-r}(x) \right)^{i_1}
\left( \frac{1}{6} f_{l,-r}(x) \right)^{i_2}
\left( \frac{1}{6} f_{l+1,-r}(x) \right)^{i_3}\\
&=&\frac{1}{2}\Big(1-\frac{f_{l-1,-r}(x)+f_{l,-r}(x)+f_{l+1,-r}(x)}{6}\Big)^{-1}\, .
\end{eqnarray*}
From this identity, we get the following recurrence relation for $l> -r$:
\begin{equation*}
2 f_{l, -r} (x) = 1 + \frac{1}{12} 2 f_{l,-r} (x) \left( 2 f_{l-1,-r}
(x) + 2 f_{l,-r} (x) + 2 f_{l+1,-r} (x) \right).
\end{equation*}
To solve this equation we follow \cite{BDFG04}. Putting
\[
F(x,y) = xy \left( 1 - \frac{1}{12} x - \frac{1}{12} y \right) -x - y,
\]
then $F \left(2 f_{l,-r} (x), 2 f_{l+1,-r} (x) 
\right)$ does not depend on $x \in [0,1]$ and $l \geq - r$, since
\begin{align*}
F & \left(2 f_{l,-r} (x), 2 f_{l+1,-r} (x)  \right) - F \left(2
f_{l-1,-r} (x), 2 f_{l,-r} (x)  \right) = \left(2 f_{l+1,-r} (x) -
2 f_{l-1,-r} (x) \right) \\
& \qquad \qquad \times
\left(2 f_{l,-r} (x) - 1 - \frac{1}{12} 2 f_{l,-r} (x) \left( 2
f_{l-1,-r} (x) + 2 f_{l,-r} (x) + 2 f_{l+1,-r} (x) \right) \right).
\end{align*}
It is easy to verify that $f_{l,-r} (x) \to 1$ as $l \to \infty$  and
since $F(2,2) = -\frac{4}{3}$, we have the relation:
\begin{equation}
\label{recgen}
F \left( 2 f_{l,-r} (x), 2 f_{l+1,-r} (x) \right) = - \frac{4}{3},
\end{equation}
for $l \geq -r$, with the initial condition $f_{-r,-r} (x) = x$. The general
solution of \eqref{recgen} is given by
\begin{equation}
\label{gen}
f_{l,-r} = 1 - \frac{2}{(l + r + a) (l + r + 1 +a)},
\end{equation}
for $l \geq -r$, where $a=a(x)$ is a function, which from the initial condition is found to be
\[
a(x) = \frac{-1 + \sqrt{1 +8 \left(1 - x \right)^{-1}}}{2},
\]
for $x \in [0,1]$.

\bigskip

Substituting \eqref{gen} in \eqref{laplaceexit}, one gets:
 \begin{eqnarray}
 &&\exp \left( \frac{2\lambda}{r^2} \right)
E  \left[ \exp - 2\lambda \frac{|\partial F_{-r}(Q_\infty)|}{r^2}\right] \nonumber \\
 &&=  E \left[ \exp 
2 \sum_{i= 0}^{\sigma_{r}-1}
\log \left(1 - \frac{2}{\left(X_i + r  + a(e^{ -
2\lambda/r^2}) \right) \left( X_i + r +1+ a(e^{ -
2\lambda/r^2})\right)} \right) 
\right]\nonumber \\
 &&= E \left[ \exp 2r^2
\int_0^{ \sigma_{r} / r^2}
\mathrm{d}t 
 \log \left(1 -
\frac{2}{\left(X_{\lfloor r^2 t \rfloor} + r + a(e^{ -
2\lambda/r^2})\right) \left( X_{\lfloor r^2 t \rfloor} + r +1+
a(e^{ -
2\lambda/r^2})\right)} \right)
\right]. \nonumber \\
\label{aremaplcer}
 \end{eqnarray}
By Skorokhod's representation theorem we can find a sequence of processes $((X_{k}^{(r)})_{k\geq0},r\geq 0)$ such that for each $r\geq 1$ we have $(X_{k}^{(r)})_{k\geq 0} = (X_{k})_{k \geq 0}$ in law and such that we have the following almost sure convergence 
$$\left(\frac{1}{r}\sqrt{ \frac{3}{2}} X_{\lfloor r^2t\rfloor}^{(r)}\right)_{t \geq0} \quad \xrightarrow[r\to\infty]{a.s.}\quad (\beta_{t})_{t \geq 0},$$ where $\beta$ is a standard Brownian motion. It is also easy to check that, as $r \to
\infty$, one has $r^{-1} a(e^{ - 2\lambda/r^2}) \to
1/\sqrt{\lambda}$. This gives
\begin{align*}
r^2 & \log \left(1 - \frac{2}{\left( X^{(r)}_{\lfloor r^2 t \rfloor}
+ r  + a(e^{ -  2\lambda/r^2 }) \right) \left(
X^{(r)}_{\lfloor r^2 t \rfloor} + r +1+ a(e^{ -
2\lambda/r^2})\right)} \right) \\
& = r^2 \log \left(1 - \frac{1}{r^2} \frac{2}{\left(X^{(r)}_{\lfloor
r^2 t \rfloor}/r + 1 + a( e^{ - 2\lambda/r^2} ) /r
\right)^2} + o \left( \frac{1}{r^2}\right) \right)\\
&\longrightarrow
- \frac{2}{\left( \sqrt{\frac{2}{3}} \beta_t
+ 1 + 1/\sqrt{\lambda} \right)^2}
\end{align*}
almost surely as $r \to \infty$. Furthermore, if we denote $\sigma_{r}^{(r)}$ the first hitting time of $-r$ by the process $X^{(r)}$, then one has the almost sure convergence:  \begin{eqnarray*} \frac{\sigma_{r}^{(r)}}{r^2} &\xrightarrow[r\to\infty]{}& T_{\sqrt{3/2}}(B),  \end{eqnarray*}
where $T_{\sqrt{3/2}}(B)$ is the first hitting time of $-\sqrt{3/2}$ of the Brownian motion $B$. This is an easy consequence of the fact that almost-surely, $B$ takes values strictly less than $-3/2$ on any time-interval of the form $[T_{\sqrt{3/2}}(B),T_{\sqrt{3/2}}(B)+\varepsilon]$ for $\varepsilon>0$.  Thus replacing $X$ by $X^{(r)}$ into \eqref{aremaplcer}, an argument of dominated convergence then gives 

  \begin{eqnarray*}
 E  \left[ \exp - \lambda \frac{2|\partial F_{-r}(Q_\infty)|}{r^2}\right]
 & \xrightarrow[r\to\infty]{} & E \left[ \exp \left( -4 \int_0^{T_{\sqrt{3/2}}(B)} \frac{ \mathrm{d}u}{ \left(\sqrt{\frac{2}{3}} \beta_u +1+ \sqrt{1/\lambda}\right)^2} \right) \right],  \end{eqnarray*}
and the scaling property of the Brownian motion shows that the right hand side of the last display is equal to

$$E \left[\exp - 6 \int_0^{T_{1}(B)} \frac{\mathrm{d}t}{\left( B_t + 1 +
\sqrt{1/\lambda} \right)^2} \right].
$$
Let us write $\mathbb{P}_x$ for the law of $(x+B_t,t\geq 0)$, and let $(Z_t,t\geq 0)$ be the canonical process. Let also $T_y(Z)=\inf\{t\geq 0:Z_t=y\}\in [0,\infty]$ be the first hitting time of $y$.  By translation, we can re-write the previous expectation as
$$\E_{1+\sqrt{1/\lambda}}\Big[\exp\Big(-6\int_0^{T_{\sqrt{1/\lambda}}(Z)}\frac{\d t}{Z_t^2}\Big)\Big]\, .$$
At this point, we can use the absolute continuity relations between Bessel
processes with different indices, due to Yor \cite[Exercise
XI.1.22]{RY99} (see also \cite{LGW06} for a similar use of these
absolute continuity relations).  The last expectation then equals
$$\Big(\frac{1+\sqrt{1/\lambda}}{\sqrt{1/\lambda}}\Big)^4
\mathbb{P}_{1+\sqrt{1/\lambda}}^{(9)}(T_{\sqrt{1/\lambda}}<\infty)\, ,$$ where $\mathbb{P}_x^{(9)}$ is the law of the $9$-dimensional Bessel process started from $x>0$.  It is classical that $\mathbb{P}_x^{(9)}(T_y<\infty)=(y/x)^{7}$ for every positive $x,y$ with $y<x$. This can be verified from the fact that $(Z_t^{-7},t\geq 0)$ is a local martingale under $\mathbb{P}_x^{(9)}$, as can be checked from Ito's formula, and the fact that $Z$ under $\mathbb{P}_x^{(9)}$ has same distribution as the Euclidean norm of a $9$-dimensional Brownian motion started from a point with norm $x$. We finally obtain that 
$$E \left[\exp - 6 \int_0^{T_{1}(B)} \frac{\mathrm{d}t}{\left( \beta_t + 1 +
      \sqrt{1/\lambda} \right)^2} \right]=\Big(\frac{1}{1+\sqrt{\lambda}}\Big)^3\, ,$$ as wanted.  \endproof

\section{Random walk on the UIPQ} \label{SRW}
This section focuses on the simple random walk over the UIPQ. We first provide a proof of a known fact (see \cite{Kri08}) that the distribution of the UIPQ is invariant under re-rooting along a simple random walk. We then make a step in understanding the recurrence/transience property of the walk on the UIPQ.

\subsection{Invariance under re-rooting along the random
  walk}\label{sec:invariance-under-re}

Let $q$ be a rooted quadrangulation, which can be finite or
infinite. We consider the nearest-neighbor random walk on $q$
starting from $e^*_+$. Rather than the random sequence of vertices
visited by this walk, we really want to emphasize the sequence of
edges that are visited. Formally, we consider a random infinite
sequence of oriented edges $(E_0,E_1,E_2,\ldots)$ starting with the
root edge $E_0=e^*$ and defined recursively as follows. Conditionally
given $(E_i,0\leq i\leq j)$, we let $E_{j+1}$ be a random edge
pointing from $(E_j)_+$, chosen uniformly among the $\deg( (E_j)_+)$
possible ones. The sequence $((E_{1+i})_-,i\geq 0)$ is then the usual
nearest-neighbor random walk on $V(q)$, starting from $e^*_+$.

We let $P_q$ be the law of the sequence
$(E_i,i\geq 0)$\footnote{Recall that a map is an equivalence class of embedded graphs, so the last definition does not really make sense but the reader can check that all quantities computed in the sequel do not depend on a representative embedded graph of the map.}. Also, for any oriented edge $e$ of the map $q$, we let $q^{(e)}$ be the map $q$ re-rooted at $e$. Finally, if $\lambda$ is a probability distribution on $\carte$, let $\mathsf{\Theta}^{(r)}(\lambda)$ be the probability distribution defined by 
$$\mathsf{\Theta}^{(r)}(\lambda)(A) = \int_{\carte}\lambda(\d q)\int P_{q}(\op{d}(e_{0},e_{1},e_{2},\ldots)) \ind_{q^{(e_{r})} \in A}, $$ for any Borel subset $A$ of $\carte$. The probability measure $\mathsf{\Theta}^{(r)}(\lambda)$ is the distribution of a random map with distribution $\lambda$, re-rooted at the $r$th step of the random walk.

\begin{proposition}
  \label{unimodular} The law $\nu$ of the UIPQ is invariant under
  re-rooting along a simple random walk, in the sense that for every
  $r\geq 0$, one has $\mathsf{\Theta}^{(r)}(\nu)=\nu$.

  Moreover, if $A$ is an event of the Borel $\sigma$-algebra of
  $(\carte,d_{\carte})$ such that $\nu(A)=1$, then
$$\nu\left(\left\{q\in \carte:\forall\, e\in 
    \overrightarrow{E}(q), \,  q^{(e)}\in A\right\}\right)= 1 \, .$$
\end{proposition}

See \cite{AL07,BCstationary} for a general study of random graphs that are invariant under re-rooting along the simple random walk. In the case of the UIPQ, the first assertion of Proposition \ref{unimodular} appears in \cite[Section 1.3]{Kri08}, see also \cite[Theorem 3.2]{AS03} for a similar result in the case of the
UIPT. We provide a detailed proof for the sake of completeness. 

\proof It is easy to see that the function $\mathsf{\Theta}^{(r)}$ on the set
$\mathcal{P}(\carte)$ of Borel probability measures on
$(\carte,d_{\carte})$ coincides with the $r$-fold composition of
$\mathsf{\Theta}=\mathsf{\Theta}^{(1)}$ with itself. Therefore, it suffices to show the result for
$r=1$. 

Let us check that $\mathsf{\Theta}$ is continuous when $\mathcal{P}(\carte)$ is
endowed with the topology of weak convergence. Indeed, if $\lambda_n$
converges weakly to $\lambda$ as $n\to\infty$, then by the Skorokhod
representation theorem, we can find a sequence $(Q_n,n\geq 0)$ of
random variables in $\carte$ with respective laws $(\lambda_n,n\geq
0)$, that converges a.s.\ to a random variable $Q$ with law
$\lambda$. For every fixed $R>0$, it then holds that
$B_{\carte,R}(Q_n)=B_{\carte,R}(Q)$ for every $n$ large
enough a.s.. Now, we can couple in an obvious way the random walks with
laws $P_{Q_n}$ and $P_Q$, in such a way that the first step $E_1$ is
the same edge in $Q_n$ and $Q$ on the event where
$B_{\carte,1}(Q_n)=B_{\carte,1}(Q)$. For such a coupling, we then
obtain that $B_{\carte,R-1}(Q_n^{(E_1)})=B_{\carte,R-1}(Q^{(E_1)})$
for every $n$ large enough. Since $R$ is arbitrary, this shows that
$Q_n^{(E_1)}$ converges a.s.\ to $Q^{(E_1)}$, so that
$\mathsf{\Theta}(\lambda_n)$ converges weakly to $\mathsf{\Theta}(\lambda)$, as desired.

Since we know by Theorem \ref{krikun} that the uniform law $\nu_n$ on
$\carte_n$ converges to $\nu$, it suffices to show that
$\mathsf{\Theta}(\nu_n)=\nu_n$. Now consider the law of the doubly-rooted map
$(q,e^*,e_1)$ under the law $\nu_n(\d q)P_q(\d(e_i)_{i\geq 0})$. The
probability that $(q,e^*,e_1)$ equals a particular doubly-rooted map
$(q,e',e'')$ with $e'_+=e''_-$ is equal to $(\#\carte_n
\deg(e'_+))^{-1}$, from which it immediately follows that
$(q,e_*,e_1)$ has the same distribution as
$(q,\overleftarrow{e}_1,\overleftarrow{e}_*)$, still under $\nu_n(\d
q)P_q(\d(e_i)_{i\geq 0})$. Hence $(q,\overleftarrow{e}_1)$ under
$\nu_n(\d q)P_q(\d(e_i)_{i\geq 0})$ has the same law $\nu_n$ as $(q,e_*)$. Since $\nu_n$ is obviously invariant under the reversal of the
root edge, we get that $(q,e_1)$ has law $\nu_n$. But by definition,
it also has law $\mathsf{\Theta}(\nu_n)$, which gives the first assertion of
Proposition \ref{unimodular}. 

Let us now prove the last part of the statement of the proposition. By
the first part, we have 
\begin{eqnarray*}
 \int_\carte \nu( \mathrm{d}q) E_q\left[\sum_{n=0}^\infty \ind_{A^c}(q^{(e_n)})\right] & =& 0.
 \end{eqnarray*}
 Thus, $\nu( \mathrm{d}q)$ a.s.\ , $E_q[\sum_{n=0}^\infty \ind_{A^c}(q^{(e_n)})] =0$. But
 \begin{eqnarray*}
E_q\left[\sum_{n=0}^\infty \ind_{A^c}(q^{(e_n)})\right] & \geq & \sum_{e\in \overrightarrow{E}(q)} P_q(\exists n \geq 0 : e_n = e) \ind_{A^c}(q^{(e)}), \end{eqnarray*} and $P_q(\exists n \geq 0 : e_n = e) >0$ for every $e\in \overrightarrow{E}(q)$ because $q$ is connected. This completes the proof.\endproof

\begin{remark}It can seem a little unnatural to fix the first step of the random
walk to be equal to $e^*$, hence to be determined by the rooted map
$q$ rather than by some external source of randomness. In fact, we
could also first re-root the map at some uniformly chosen random edge
incident to $e^*_-$, and start the random walk with this new
edge. Since the first re-rooting leaves the laws $\nu_n,\nu$
invariant, as is easily checked along the same lines as the previous
proof, the  results of Proposition \ref{unimodular} still hold with the new random walk.
\end{remark}
\subsection{On recurrence}

Let $Q_{\infty}$ be the uniform infinite planar quadrangulation. Conditionally on $Q_{\infty}$,  $(E_{k})_{k\geq0}$ denotes the random sequence of oriented edges with $E_{0}=e^*$ traversed by a simple random walk on $Q_{\infty}$ as discussed at the beginning of Section \ref{sec:invariance-under-re}. We write $ \mathcal{X}_{k}=(E_{k})_{-}$ for the sequence of vertices visited along the walk. For $k\geq 0$, we denote  the quadrangulation $Q_{\infty}$ re-rooted at the oriented edge $E_{k}$ by $Q_{\infty}^{(k)}$.  Proposition \ref{unimodular} shows that $Q_{\infty}^{(k)}$ has the same distribution as $Q_{\infty}$. 

\begin{open}[\cite{AS03}] Is the simple random walk $(\mathcal{X}_{k})_{k\geq0}$ on $Q_{\infty}$ almost surely recurrent?
\end{open}
A similar question for UIPT arose when Angel \& Schramm \cite{AS03} introduced this infinite random graph. These questions are still open. James T. Gill and Steffen Rohde \cite{RG10} proved that the Riemann surface obtained from the UIPQ by gluing squares along edges is recurrent for Brownian motion. The first author and Itai Benjamini also proved that the UIPQ is almost surely Liouville \cite{BCstationary}. However the lack of a bounded degree property for the UIPQ prevents one from deducing recurrence from these results (see also \cite{BS01}). Our new construction of the UIPQ however leads to some new information suggesting that the answer to the above Question should be positive.

\begin{theorem} \label{rec} The process $(\ell(\mathcal{X}_{n}))_{n\geq 0}$ is a.s.\ recurrent, $i.e.$\ visits every integer infinitely often. \end{theorem}

\proof For every $k\geq0$, one can consider the
labeling $(\ell^{(k)}(u))_{u\in Q_\infty}$ of the vertices of
$Q_\infty$ that corresponds to the labeling given by Theorem \ref{sec:unif-infin-label-1} applied to the rooted infinite planar quadrangulation $Q_\infty^{(k)}$. On the one hand, it is straightforward to see from \eqref{eq:4} that $\ell^{(k)}(u)-\ell^{(k)}(v) =\ell(u)-\ell(v)$ for every $u,v \in Q_\infty$. On the other hand, applying Proposition \ref{unimodular} we deduce that  the process $(\ell^{(k)}(\mathcal{X}_{k+i})-\ell^{(k)}(\mathcal{X}_k))_{i\geq 0}$ has the same distribution as $(\ell(\mathcal{X}_i)-\ell(\mathcal{X}_0))_{i\geq0}.$ Gathering up the pieces, we deduce that for every integer $k \geq 0$ we
have
\begin{eqnarray}
\label{stationary}
\big(\ell(\mathcal{X}_{i})-\ell(\mathcal{X}_0)\big)_{i\geq 0} &\overset{(d)}{=}&
  \big(\ell(\mathcal{X}_{k+i})-\ell(\mathcal{X}_{k})\big)_{i\geq 0}.
\end{eqnarray}
Hence the increments $(\ell(\mathcal{X}_{i+1})- \ell(\mathcal{X}_{i}))_{i\geq 0}$ form is a stationary sequence.  Furthermore, we have $|\ell(\mathcal{X}_{1})-\ell(\mathcal{X}_{0})| =1,$ and since the distribution of $Q_{\infty}$ is preserved when reversing the orientation of the root edge we deduce 
 $$ \ell(\mathcal{X}_{1})-\ell(\mathcal{X}_{0}) \quad \overset{(d)}{=}\quad \ell(\mathcal{X}_{0})-\ell(\mathcal{X}_{1}) \quad\overset{(d)}{=} \quad\mathcal{B}(1/2).$$ In particular the increments of $\ell(\mathcal{X}_{n})$ have zero mean. Suppose for an instant that the increments of $\ell(\mathcal{X}_{n})$ were also ergodic, then Theorem 3 of \cite{Dek82} would directly apply and give the recurrence of $\ell(\mathcal{X}_{n})$. Although the UIPQ is ergodic, a proof of this fact would take us too far, so we will reduce the problem to the study of ergodic components.

 By standard facts of ergodic theory, the law $\xi$ of the sequence of increments $(\ell(\mathcal{X}_{i+1})-\ell(\mathcal{X}_{i}))_{i\geq0}$ can be expressed as a barycenter of ergodic probability measures in the sense of Choquet, namely for every $A\subset \mathcal{B}(\mathbb{R})^{\otimes \mathbb{N}}$ we have
\begin{eqnarray} \xi(A) &=& \int  \zeta(A) \mathrm{d}m(\zeta), \label{choquet}
\end{eqnarray}
where $m$ is a probability measure on the set of all probability measures on $(\mathbb{R}^\mathbb{N}, \mathcal{B}(\mathbb{R})^{\otimes \mathbb{N}},\mathbb{P})$ that are ergodic for the shift. In our case, it suffices to show that $m$-almost every $\zeta$ satisfies the assumption of \cite[Theorem 3]{Dek82}. Specializing \eqref{choquet} with $A_{1} =  \{ (y_{i})_{i\geq 0} : |y_{i+1}-y_{i}| \leq 1, \forall i \geq 0\}$ we deduce that $m$-almost every $\zeta$, we have $\zeta(A_{1})=1$, in particular the increments under $\zeta$ are integrable. It remains to show that they have zero mean. 
\begin{lemma} \label{Carne} Almost surely we have 
\begin{eqnarray*}   \lim_{n\to \infty} \frac{\ell(\mathcal{X}_{n})}{n} & =& 0.
\end{eqnarray*}
\end{lemma}
 \proof In \cite[Theorem 6.4]{CD06} it is shown that $\mathbb{E}[\#B_{\carte,r}(Q_{\infty})] \leq C_3 r^4$ where $C_3>0$ is independent of $r\geq 1$. Using the  Borel-Cantelli lemma we easily deduce that  \begin{eqnarray} \label{upball} \lim_{r\to \infty} r^{-6}\#B_{\carte,r}(Q_{\infty}) &=&0, \quad  a.s. \end{eqnarray}
We now use the classical Varopoulos-Carne upper bound (see for instance Theorem 13.4 in \cite{LP10}): we have 
 \begin{eqnarray} p_{n}(e^*_{+},x) &\leq& 2 \sqrt{\frac{\op{deg}(x)}{\op{deg}(e^*_{+})}} \exp\left( -\frac{\op{d}_{\op{gr}}^{Q_{\infty}}(e^*_{+},x)^2}{2n} \right), \label{varo-carne} \end{eqnarray}where conditionally on $Q_{\infty}$, $p_{n}(.,.)$ is the $n$-step transition probability of the simple random walk started from $e^*_{+}$ in $Q_{\infty}$.  Conditionally on $Q_\infty$, using a crude bound $\op{deg}(x) \leq \# B_{\carte,n+1}(Q_{\infty})$ on the degree of a vertex $x \in B_{\carte,n}(Q_{\infty})$,  we have using \eqref{varo-carne}  
\begin{eqnarray*} \label{upkernel} P_{Q_\infty}(\mathcal{X}_{n} \notin  B_{\carte,n^{2/3}} (Q_{\infty})) &\leq& 2 \exp\left(-\frac{n^{1/3}}{2}\right) \big(\#B_{\carte,n+1}(Q_\infty)\big)^{3/2}.\end{eqnarray*} Hence on the event $\{\lim_{r\to \infty} r^{-6}\#B_{\carte,r}(Q_{\infty}) =0\}$, an easy application of   the  Borel-Cantelli lemma shows that $n^{-1} \op{d}_{\op{gr}}^{Q_{\infty}}(\mathcal{X}_{n},\varnothing) \to 0$ as $n \to \infty$. Since $|\ell(\mathcal{X}_{n})| \leq \op{d}_{\op{gr}}^{Q_{\infty}}(\mathcal{X}_{n},\varnothing)$, the above discussion together with \eqref{upball}  completes the proof of the lemma.\endproof

\noindent Let us complete the proof of Theorem \ref{rec}. We can specialize formula  (\ref{choquet}) to $A_{2} = \{(y_{i})_{i\geq 0} : \lim i^{-1} |y_{i}|=0 \},$ to obtain that $m$-a.e $\zeta$ we have $\zeta(A_{2})=1$. Using the ergodic theorem that means that the increments under $\zeta$ are centered. We can thus apply Theorem 3 of \cite{Dek82} to get that for $m$-almost every $\zeta$, the process whose increments are distributed according to $\zeta$ is recurrent, hence $(\ell(\mathcal{X}_n))$ is almost surely recurrent. \endproof

\section*{Appendix: infinite maps and their embeddings}

In this section, we explain how the elements of $\carte_\infty$ can be seen as infinite quadrangulations of a certain non-compact surface, completing the description of Section \ref{sec:infin-quadr}. 

Recall that an element $q$ of $\mathbf{Q}_\infty$ is a sequence of compatible maps with holes $(q_1,q_2,\ldots)$, in the sense that $q_r=B_{\carte,r}(q_{r+1})$. This sequence defines a unique cell complex $S_q$ up to homeomorphism, with an infinite number of 2-cells, which are quadrangles. This cell complex is an orientable, connected, separable topological surface, and every compact connected sub-surface is planar. 

It is known \cite{Ric63} that the topology of $S_q$ is characterized by its {\em ends space}, which is a certain totally disconnected compact  space. Roughly speaking, the ends space determines the different ``points at infinity'' of the surface. More precisely, following \cite{Ric63}, we define a boundary component of $S_q$ as a sequence $(U_1,U_2,\ldots)$ of subsets of $S_q$, such that 
\begin{itemize}
\item for every $i\geq 1$, the set $U_i$ is unbounded, open, connected and with compact boundary, 
\item for every $i\geq 1$, it holds that $U_{i+1}\subset U_i$, 
\item for every bounded subset $A\subset S_q$, $U_i\cap A=\varnothing$ for every $i$ large enough. 
\end{itemize}
Two boundary components $(U_i,i\geq 1),(U'_i,i\geq 1)$ are called equivalent if for every $i\geq 1$ there exists $i'\geq 1$ such that $U'_{i'}\subset U_i$, and vice-versa. An {\em end} is an equivalence class of boundary components. For every $U\subset S_q$ with compact boundary, we let $V_U$ be the set of all ends whose corresponding boundary components are sequences of sets which are eventually included in $U$. The topological space having the sets $V_U$ as a basis is called the ends space, and denoted by $\mathscr{E}_q$.  

Conversely, it is plain that every rooted quadrangulation of an orientable, connected, separable, non-compact planar surface, defines an element of $\carte_\infty$, by taking the sequence of the balls centered at the root vertex, with the same definition as in Section \ref{sec:infin-quadr}. The separability ensures that the collection of balls exhausts the whole surface. Thus we have:

\begin{proposition}
  The elements of $\carte_\infty$ are exactly the quadrangulations of orientable, connected, separable, non-compact planar surfaces, and considered up to homeomorphisms that preserve the orientation.  
\end{proposition}

To understand better what the ends space is in our context, note that there is a natural tree structure $\mathscr{T}_q$ associated with $q\in \carte_\infty$. The vertices $v$ of this tree are the holes of $q_1,q_2,q_3,\ldots$, and an edge links the vertices $v$ and $v'$ if there exists an $r\geq 1$ such that $v$ is a hole of $q_r$, $v'$ is a hole of $q_{r+1}$, and $v'$ is included in the face determined by $v$. Furthermore, all the holes in $q_1$ are linked by an edge to an extra root vertex. 

It is then easy to see that $\mathscr{E}_q$ is homeomorphic to the ends space $\partial\mathscr{T}_q$ which is defined as follows: $\partial\mathscr{T}_q$ is just the set of infinite injective paths (spines) in $\mathscr{T}_q$ starting from the root, and a basis for its topology is given by the sets $W_v$ made of the spines that pass through the vertex $v$ of $\mathscr{T}_q$. 
(This is consistent, since it is easy and well-known that the ends space of trees with finite degrees is a compact totally disconnected space.)

In particular, when $\mathscr{T}_q$ has a unique spine, then $\mathscr{E}_q$ is reduced to a point, which means that the topology of $S_q$ is that of the plane $\mathbb{R}^2$.

\addcontentsline{toc}{section}{References}
\bibliographystyle{abbrv}
\bibliography{bibli}
\end{document}